\documentclass[12pt,a4paper]{article}
\makeatletter \@addtoreset{equation}{section} \makeatother
\def\N{\mbox{I\hspace{-.15em}N}}
\def\R{\mbox{I\hspace{-.15em}R}}

\newtheorem{theo}{Theorem}[section]
\newtheorem{lem}[theo]{Lemma}

\newtheorem{pro}[theo]{Proposition}
\def\u{\mathbf{u}}
\def\U{\mathbf{U}}
\def\o{\mbox{\boldmath$\omega$\unboldmath}}

\def\vfy{\varphi}

\def\tho{\mbox{\boldmath$\tau$\unboldmath}}
\def\un{\underline}
\def\O{\Omega}
\def\v{\mathbf{v}}

\def\p{\Vert}

\def\n{\nabla}
\def\d{\mathrm{div}\,}

\def\P{\partial}
\def\F{\frac}

\def\a{\alpha}
\def\m{\mb{m}}

\def\t{\times}

\def\W{\mathcal{W}}
\def\z{\mathbf{z}}

\def\hs{\hspace*{0.2cm}}
\def\mb{\mathbf}
\def\mr{\mathrm}
\def\H{\hspace*{0.7cm}} 
\def\iy{\infty}

\def\di{\displaystyle}
\def\rk{$\rbrack$}
\def\lk{$\lbrack$}

\def\vb{\mbox{\boldmath$\vfy$\unboldmath}}

\def\tb{\textbf}
\def\G{\Gamma}
\def\x{\mb{x}}
\def\ov{\overline}
\def\y{\mb{y}}
\def\tho{\mbox{\boldmath$\tau$\unboldmath}}
\usepackage{color}
\usepackage{graphics}
\usepackage{graphicx}
\begin{document}
\setcounter{page}{0}
\title{Solutions in $H^1$ of the steady transport equation in a \mbox{bounded polygon
with a fully non-homogeneous velocity.}}
\author{J. M. Bernard*}
\date{ }
\maketitle
\begin{abstract} This article studies the solutions in $H^1$ of a
  steady transport equation with a divergence-free driving velocity
  that is $W^{1,\iy}$, in a two-dimensional bounded polygon. Since the velocity is assumed
  fully non-homogeneous on the boundary, existence and uniqueness of
  the solution require a boundary condition on the open part $\Gamma^-$ where the normal component of $\u$ is strictly 
negative. In a previous article, we studied the solutions in $L^2$ of this steady transport equation. The methods, developed in this article, can be extended to prove existence and uniqueness of a solution in $H^1$ with Dirichlet boundary condition on $\Gamma^-$ only in the case where the normal component of $\u$ does not vanish at the boundary of $\Gamma^-$. In the case 
where the normal component of $\u$ vanishes at the boundary of $\Gamma^-$, under appropriate assumptions, we construct  local $H^1$ solutions in the neighborhood of the end-points of $\Gamma^-$, which allow us to establish existence and uniqueness of the solution in $H^1$ for the transport equation with a Dirichlet boundary condition on $\Gamma^-$.
\end{abstract}  
 \vspace{1cm}
\renewcommand{\abstractname}{R\'esum\'e}
\begin{abstract} Cet article \'etudie les solutions dans $H^1$ d'une \'equation de transport stationnaire avec une vitesse de r\'egularit\'e 
$W^{1,\iy}$ \`a divergence nulle, dans un polygone born\'e. La vitesse \'etant suppos\'ee non nulle sur la fronti\`ere, l'existence et l'unicit\'e de la solution requi\`erent une condition sur la partie de la fronti\`ere o\`u la composante normale de la vitesse est strictement n\'egative. Dans un pr\'ec\'edent article, nous avons \'etudi\'e les solutions dans $L^2$ de cette \'equation de transport stationnaire. Les m\'ethodes, d\'evelopp\'ees dans cet article, peuvent \^{e}tre \'etendues pour prouver l'existence et l'unicit\'e d'une solution dans $H^1$ avec une condition de Dirichlet sur $\G^-$ seulement dans le cas o\`u la composante normale de $\u$ ne s'annulle pas \`a la fronti\`ere de $\G^-$. Dans le cas o\`u la composante normale de $\u$ s'annulle \`a la fronti\`ere de $\G^-$, sous des hypoth\`eses appropri\'ees, nous construisons des solutions locales au voisinage des points fronti\`eres de $\G^-$ de r\'egularit\'e $H^1$, qui nous permettent d'\'etablir l'existence et l'unicit\'e de la solution dans $H^1$ de l'\'equation de transport avec une condition de Dirichlet sur $\G^-$.
\end{abstract}
\vspace{0.3cm}
\tb{Key words.} transport equation, nonstandard boundary condition, localization methods\\[0.3cm]
*Laboratoire Analyse et Probabitilit\'es, Universit\'e d'Evry Val
d'Essonne,\\
 23 Boulevard de France, 91037 Evry, France.\\
e-mail: jm-bernard@club-internet.fr
  \pagebreak
  \setcounter{section}{-1.}
  \section{Introduction.}
  \H Transport equations are studied in many frameworks. In [2,7,12] the stress $z$, i.e, the transported quantity, is not assumed regular, while they impose strong conditions on the fluid velocity $\u$, which indicates the direction of the transport. Contrary to this, in [1,4], the velocity has only bounded variation 
  with its divergence integrable, but the stress is assumed bounded or continuous. In fact, we have to choose the regularity of $z$ and $u$ for the product $\u\,.\,\n z$ to be well defined in some distributional sense. Thus, in [7], V. Girault and L.R. Scott, for defining $\u\,.\,\n z$ with the weaker assumptions, studied the transport equation with the stress $z$ in $L^2(\O)$, the velocity $\u$ given in $H^1(\O)^d$, with $\d\u=0$, and the right hand side given in $L^2(\O)$, where $\O$ is a Lipschitz-continuous domain. These Authors established existence and uniqueness of the solution for the transport equation by using the essential technique of Puel and Roptin [11] and the renormalizing argument of DiPerna and Lions [5]. In a following article [8], they extended their results from 
  $L^2$ to $H^1$ for the transport equation. By another technique, in particular a Yosida aproximation, V. Girault anf L. Tartar [9] studied the solutions in $L^p$, $p\ge 2$, of the transport equation, when the right hand side is in $L^p$.\\
\H  However, all these approaches of transport equations assume that the normal component of the fluid velocity $\u$ vanishes on the boundary of the domain. Indeed, in the contrary case, the problem is no longer well-posed and the unicity requires a boundary condition. However, it is not possible to define the trace on the boundary of the stress $z$ when it is not regular but only square-integrable. Nevertheless, such transport equation with $\u\,.\,\mb{n}\not= 0$, where $\mb{n}$ denotes the unit exterior normal to the boundary, arises in the problem of fully nonhomogeneous second grade fluid [7]: multiple solutions imply that additional boundary conditions should be imposed.\\
  \H In a previous article [3], we established existence and uniqueness of the solution, in the space where $z$ and $\u\,.\,\n z$ are $L^2$, for the transport equation, with a boundary condition on the open part of the boundary where the normal component of $\u$ is strictly negative, where  $\O$ is a Lipschitz-continuous domain of $\R^d$, $\u$ is given in $H^1(\O)^d$ such that $\d\u=0$, the right hand side is given in $L^2(\O)$, and $\mathcal{W}$ is a given real parameter different from 0. We showed that it is possible to define the normal component of $z\u$ on the boundary and, hence, to prove that the problem is well-posed by requiring a condition for the normal component of $z\u$ on the part of the boundary where $\u\,.\,\mb{n}<0$.\\
 \H The present article studies the steady transport problem : find $z\in H^1(\O)$ such that
  \begin{equation}\label{1pbtH1}\left\lbrace \begin{array}{l}z+\mathcal{W}\u\,.\,\n z=l\H \mr{in}\ \O,\\
z=0\hspace*{2.7cm} \mr{on}\ \G^- \end{array}\right.,\end{equation}
where  $\O$ is a bounded polygon of $\R^2$, $\u$ is given in $W^{1,\iy}(\O)^2$ such that $\d\u=0$, $\Gamma^-$ is the open part of the boundary of $\O$ such that $\u\,.\,\mb{n}<0$, $l$ is given in $H^1(\O)$, and $\mathcal{W}$ is a given real parameter different from 0. But, in a such framework, if we look for  a solution in $H^1$, a difficulty arises when $\u\,.\,\mb{n}$ vanishes at the boundary of $\Gamma^-$, as we shall see in examples given below. Indeed, 
 the fact that the function $\u.\mb{n}$ vanishes at a boundary point of $\G^-$ leads to a discontinuity for the 
partial derivatives of the solution $z$ at this point and the solution $z$ has not always the regularity $H^1$, see examples 4, 5 and 7. As we will see in the following examples, the regularity of the solution seems to depend on the multiplicity of the root of the equation $\u\,.\,\mb{n}=0$ at the boundary of $\G^-$ and on the sign of $\u\,.\,\tho_{-}$, where $\tho_{-}$ is the unit tangent vector to the boundary at the point $m$, directed towards $\G^-$ : in these examples, the solution of the transport equation is $H^1$ if the multiplicity of the root is 0 or 1 and if the sign of $\u\,.\,\tho_{-}$ is negative in $\m$. On the contrary, the solution is not $H^1$ if the multiplicity of the root is strictly greater than 1 or if the sign of $\u\,.\,\tho_{-}$ in $\m$ is positive, which is consistent with the assumptions (\ref{4ctheun}) of the Theorem \ref{5theu}. \\ 
\H When $\u\,.\,\mb{n}$ does not vanish at the boundary of $\G^-$, by using results and tools of [3], we can prove existence and uniqueness of the solution $H^1$ for the steady transport equation with the boundary condition 
on $\G^-$. In contrast, when $\u\,.\,\mb{n}$ vanishes at the boundary of $\G^-$, the previous method does not work anymore. In this case, we split the right-hand side of the transport equation, which gives us a set of localized problems, and the solution $H^1$ of the transport problem is the sum of the solutions $H^1$ of the localized problems. For solving the problems localized in the neighborhoods of the points where $\u\,.\,\mb{n}$ vanishes such as simple roots of the equation $\u\,.\,\mb{n}=0$, in the case where $\u\,.\,\tho_{-}$ is negative at each of these points, we use a change of variables, which allows us to explain the local solution $H^1$ of the transport equation in integal form. Next, we extend this local solution to the whole domain $\O$ and we obtain the $H^1$ solutions of the transport problems localized around these roots. Instead, to solve the problems localized far enough of these roots, the methods of [3] yield the $H^1$ solutions.\\[0.2cm]
\H After this introduction, this article is organized as follows. In
section 1, we study several examples of transport problems, which show the link between the regularity of the solution ($L^2$ or $H^1$) and both the multiplicity of the roots of the equation $\u\,.\,\mb{n}=0$ at the end points of $\G^-$ and the sign of $\u\,.\,\tho_{-}$ at these points. Section 2 is
devoted to the solution in $H^1$ of the transport problem when the normal component of the velocity does not vanish on $\ov{\Gamma^-}$. In section 3, we deal with the solutions in $H^1$ of the transport problem in the case where the normal component of the velocity vanishes on $\ov{\Gamma^-}$.\\[0.3cm]
    \H We end this introduction by recalling some basic results of [3] that we shall use throughout this article. Let $\Gamma'$ be an open part of the boundary $\P\O$ of class $C^{0,1}$
and, for $r>2$, $T_{1,r}^{\Gamma'}$ the mapping $v\mapsto v_{|\Gamma'}$
defined on $W^{1,r}(\O)$.
We denote by $W^{1-\F{1}{r},r}(\Gamma')$ (see [10]) the space
$T_{1,r}^{\Gamma'}(W^{1,r}(\O))$ which is equipped with the norm:
\begin{equation}\label{1dn1/2}
\p\varphi\p_{W^{1-1/r,r}(\Gamma')}=
\inf\{\p v\p_{W^{1,r}(\O)},\ v\in W^{1,r}(\O)\ \mr{and}\ v_{|\Gamma'}=\varphi\}.
\end{equation}
\H For fixed $\u$ in $H^1(\O)^2$, let us introduce the space 
 \begin{equation}\label{1dXu}
 X_{\u}(\O)=\{z\in L^2(\O),\ \u\,.\,\n z\in L^2(\O)\},\end{equation}
 which is a Hilbert space equipped with the norm
 \begin{equation}\label{1dnXu}
 \p z\p_{\u}=(\p z\p_{L^2(\O)}^2+ \p\u\,.\,\n z\p_{L^2(\O)}^2)^{1/2}.
 \end{equation}
  In the same way we define
  $$Y_{\u}(\O)=\{z\in L^2(\O),\ \u\,.\,\n z\in L^1(\O)\}.$$
  We recall a theorem ( see [3]) concerning the normal component of boundary values of $(z\u)$
where $z$ belongs to $Y_{\u}(\O)$.
\begin{theo}\label{1tg'n}
Let $\O$ be a Lipschitz-continuous domain of $\R^d$, let $\u$ belong to
$H^1(\O)^d$ with $\d\u=0$ in
$\O$ and let $r>d$ be a real number. We denote by $r'$ the
real number defined by: $\di\F{1}{r}+\di\F{1}{r'}=1$. The mapping
$\gamma'_{\mb{n}}:\ z\mapsto (z\u)\,.\,\mb{n}_{|\P\O}$ defined on
$\mathcal{D}(\overline{\O})^d$ can be extended by continuity to a linear and continuous
mapping, still denoted by $\gamma'_{\mb{n}}$, from $Y_{\u}(\O)$ into
$W^{-1/r',r'}(\P\O)$.\end{theo}
\H From this theorem and with a density argument, we derive the following
  Green's formula: let $r>d$ be a real number and let $\u$ be in
  $H^1(\O)^d$ with $\d\u=0$ in $\O$,
      \begin{equation} \label{1green}
  \forall z\in Y_{\u}(\O),\ \forall\varphi\in W^{1,r}(\O),
  \ \int_{\O}z(\u\,.\,\n\varphi)\,d\mb{x}
  +\int_{\O}\varphi(\u\,.\,\n z)\,d\mb{x}
  =<(z\u)\,.\,\mb{n},\varphi>_{\P\O}.\end{equation}
  \H Let $\G_0$ and $\G_1$ be two non empty open parts of $\P\O$ that have a finite number of connected components and verify 
  $$\G^0\cap\G^1=\emptyset,\H
\P\O=\ov{\G^0}\cup\ov{\G^{1}},$$
such that $\ov{\G^0}\cap\ov{\G^{1}}$ has a finite number of connected components.\\
\H We introduce
the space $W^{-1/r',r'}(\G_0)=(W^{1-1/r,r}_{00}(\G_0))'$, where 
\begin{equation}\label{1dW00r}
W^{1-1/r,r}_{00}(\G_0)=\{v_{|\G_0},\ v\in W^{1,r}(\O),\
v_{|\G_1}=0\},\end{equation}
and we denote $<\,.\,,\,.\,>_{\G_0}$ the duality pairing between these
two spaces. Note that if $z\in Y_{u}(\O)$, then
$(z\u)\,.\,\mb{n}_{|\G_0}\in W^{-1/r',r'}(\G_0)$ and, in the same way
as previously, we have the Green's formula : $\forall z\in Y_{\u}(\O)$, 
$\forall\vfy\in W^{1,r}(\O)$, with $\vfy_{|\G_1}=0$, 
$\forall\u\in H^1(\O)^d$ with $\d\,\u=0$ in $\O$,
\begin{equation}\label{1green2}
 \int_{\O}z(\u\,.\,\n\varphi)\,d\mb{x}
  +\int_{\O}\varphi(\u\,.\,\n z)\,d\mb{x}
  =<(z\u)\,.\,\mb{n},\varphi>_{\G_0}.\end{equation}
Then, we can define the following space : 
\begin{equation}\label{1dXug0}
X_{\u}(\G_0)=\{ z\in X_{\u},\ (z\u)\,.\,\mb{n}_{|\G_0}=0\}.\end{equation}
  \H From now on, we suppose that $d=2$ and $\O\subset\R^2$. Let us denote by $\G^-$ and $\G^{0,+}$ the following open portions of
$\P\O$
\begin{equation}\label{1dg-}
\G^-=\bigcup_{i\in I} \o_i,\end{equation}
where the sequence $(\o_i)_{i\in I}$ represents the set of the open
sets $\o_i$ of $\P\O$ such that 
$\W\,\u\,.\,\mb{n}<0$ almost
everywhere in $\o_i$.  
In the same way, 
\begin{equation}\label{1dg0+}
\G^{0,+}=\bigcup_{j\in J} \o'_j,\end{equation} 
where the open sets $\o'_j$ of $\P\O$ are such that
$\W\,\u\,.\,\mb{n}\ge 0$ almost
everywhere in $\o'_j$. 
Let us note that these definitions imply
 $$\G^-\cap\G^{0,+}=\emptyset.$$
 \H We assume that $\G^-$ and $\G^{0,+}$ have a finite number of connected components and verify
\begin{equation}\label{1dg0g1}
\P\O=\ov{\G^-}\cup\ov{\G^{0,+}},\H
\ov{\G^-}\cap\ov{\G^{0,+}}=\{\mb{m}_1,\ldots,\mb{m}_q\},
\end{equation}
where $\mb{m}_k$, $1\le k\le q$, denote points of the boundary $\P\O$.\\
\H Let us define the space $U$ by 
\begin{equation}\label{dU}
U=\{\v\in H^1(\O)^2;\ \d\v=0\}.\end{equation}
Finally, we recall  basic results of [3] that we apply in the particular case where $d=2$.
\begin{pro}\label{2pgreenfg-} Let $\O$ be a Lipschitz-continuous domain
  of $\R^2$, let $\u$ be given in $U$, defined by (\ref{dU}), and let
  $\G^-$ and $\G^{0,+}$ be defined by (\ref{1dg-}) and (\ref{1dg0+}),
  verifying (\ref{1dg0g1}). 
Let $z$ belong to
 $X_{\u}(\G^-)$ and $w$ to $X_{\u}(\G^{0,+})$ . Then, $z$ and $w$ verify the following inequalities
\begin{equation}\label{2greenfg-+} 
\int_{\O}(\mathcal{W}\u\,.\,\n z)\,z\,d\x\ge
0,\H\int_{\O}(\mathcal{W}\u\,.\,\n w)\,w\,d\x\le 0.\end{equation}\end{pro}
\H Considering the problem: for $\u$ in $U$,
$l$ in $L^2(\O)$ and $\mathcal{W}$ in $\R^*$, find $z$ in $L^2(\O)$
such that:
\begin{equation}\label{2te}
\left\lbrace \begin{array}{l}z+\mathcal{W}\u\,.\,\n z=l\H \mr{in}\ \O,\\
(z\u)\,.\,\mb{n}=0\H\H\hspace*{0.15cm} \mr{on}\ \G^- \end{array}\right..
\end{equation}
In [3], we prove the following result of existence and uniqueness in $L^2$.
\begin{theo}\label{2teeth}
Let $\O$ be a lipschitz-continuous domain of $\R^2$ and let
$\G^-$ and $\G^{0,+}$ be defined by (\ref{1dg-}) and (\ref{1dg0+}),
  verifying (\ref{1dg0g1}). For all
$\u$ in $U$, defined by (\ref{dU}), all $l$ in $L^2(\O)$ and all real numbers $\mathcal{W}$
in $\R^*$, the transport problem (\ref{2te}) has a unique 
solution $z$ in $L^2(\O)$.
\end{theo}
\section{Examples of transport problems.}
\H In this section, we study different examples of transport problems (\ref{1pbtH1}), obtained with different choices of velocities $\u$, functions $l$ as right hand side and domains $\O\subset \R\t]0,+\iy[$.\\
\H First, we choose $\u(x,y)=(x,-y)$, that verifies $\d\,\u=0$ everywhere, which cor-\linebreak responds to the following transport problem : find $z\in H^1(\O)$
satisfying
\begin{equation}\label{pbh2}\left\lbrace
\begin{array}{l} z+x\,\di\F{\P z}{\P x}-y\,\di\F{\P z}{\P y}=l\H \mr{in}\ \O,\\
z_{|\G^-}=0\end{array}\right..\end{equation}
 This problem was 
introduced in [3] in the particular case where $\O=]0,1[\t ]1,2[$. In the examples 1, 2 and 3, we study the problem (\ref{pbh2}) for different examples of functions $l$ and domains $\O$.
We can set $$\left\lbrace
\begin{array}{l}X=xy\\
Y=\ln y\end{array}\right..$$
Setting $\O_*=\{(X,Y)\in\R^2,\ (X e^{-Y},e^Y)\in \O\}$, $\G_*^-=\{(X,Y)\in\R^2,\ (X e^{-Y},e^Y)\in \G^-\}$ and $Z(X,Y)=z(x,y)$, we derive the following equivalent problem :
Find $Z\in H^1(\O_*)$ satisfying
$$\left\lbrace\begin{array}{l}
Z-\di\F{\P Z}{\P Y}=l(Xe^{-Y}, e^{Y})\\
Z_{|\G_*^-}=0\end{array}\right..$$
Hence, if $a\in [\min\limits_{(x,y)\in\ov{\O}}(y),\max\limits_{(x,y)\in\ov{\O}}(y)]$, we find the general solution of the first equation of (\ref{pbh2}):\\ 
\begin{equation}\label{2sz}z(x,y)=y\int_{\ln y}^{\ln a} e^{-t}l(xye^{-t},e^{t})\,dt+
 y C(xy),\end{equation}
where $C$ if any function in $L^2$. Thus, we have an infinity of solutions. In order to obtain a well-posed problem (see [3]), it is necessary to 
require a boundary condition on $\G^-$, which allows us to compute the function $C$.\\ \pagebreak\\
\H Second, another choice of velocity $\u$ is: $\u(x,y)=(xy^2+y,-y^2)$, as in the example 4 (or $\u(x,y)=(-xy^2-y,y^2)$, as in the example 5). Setting $X=xy^2+y$, $Y=\di\F{1}{y}$ and $Z(X,Y)=z(x,y)$, we have the equivalence
$$z+\u\,.\,\n z=l\Longleftrightarrow Z+\F{\P Z}{\P Y}=\tilde{l},$$
where $\tilde{l}(X,Y))=l(XY^2-Y,\F{1}{Y})$, which implies$$\forall (x,y)\in\O,\ z(x,y)=e^{-\F{1}{y}}\int_{\F{1}{3}}^{\F{1}{y}} e^tl((xy^2+y)t^2-t,\F{1}{t})\,dt+e^{-\F{1}{y}}C(xy^2+y),$$
where $C$ is any function in $L^2$. With the choice of $l=1$, we obtain
\begin{equation}\label{2sz2}
z(x,y)=1-e^{\F{1}{3}-\F{1}{y}}+e^{-\F{1}{y}}C(xy^2+y).\end{equation}
\H For a given velocity $\u$, we introduce the following notations : 
\begin{equation}\label{2dG0}
\G^0\ \mr{is\ the\ interior\ of\ the\ set}\ \{\x\in\P\O,\ (\u\,.\,\mb{n})(\x)=0\},\end{equation} 
\begin{equation}\label{2dG+}
\G^+\ \mr{is\ the\ interior\ of\ the\ set}\ \{\x\in\P\O,\ (\u\,.\,\mb{n})(\x)>0\}.\end{equation}
\H In the first two examples, the function $\u\,.\,\mb{n}$ does not vanish on $\ov{\G^-}$ and we verify that the regularity $H^1$ of the solution $z$ depends on the regularity $H^1$ of $l$.
\subsection{Example 1 : $\O=]0,1[\t]1,2[$,\ $l(x,y)=x^{\F{2}{3}}$, $\u(x,y)=(x,-y)$.}
\begin{minipage}[b]{10cm}
\H Let $\O$ be the square $]0,1[\t]1,2[\subset \R^2$ (see figure 1.1). 
We can verify that $$\lim_{x\to 0^+}l'_x(x,y)=+\iy,\ \mr{but}\ l\in H^1(\O).$$ 
Moreover $\G^-=\G_3$ (in red), $\G^0=\G_4$ (in green),\\ $\G^+=\G_1\cup\G_2$ (in blue).
From (\ref{2sz}) and $a=2$, in view of the boundary condition $z_{|\G^-}=0$, we derive $C=0$ and
$$\forall (x,y)\in [0,1]\t [1,2],\ z(x,y)=\F{3}{5}x^{\F{2}{3}}(1-\F{y^{\F{5}{3}}}{2^{\F{5}{3}}}).$$
We can verify that $z\in H^1(\O)$. Thus, we have taken $l\in H^1(\O)$ and we have obtained $z\in H^1(\O)$ and the problem (\ref{pbh2}) is well-posed. 
\end{minipage}
\hspace*{1cm}
\unitlength=1cm
\begin{picture}(6,6)
\put(1.13,1){\vector(0,1){6.2}}
\put(0,2){\vector(1,0){5}}
\linethickness{0.04cm}
\color{blue}
\put(1,3.9){\line(1,0){2.4}}
\put(3.4,3.9){\line(0,1){2.4}}
\color{red}
\put(3.4,6.3){\line(-1,0){2.4}}
\color{green}
\put(1,6.3){\line(0,-1){2.4}}
\color{black}
\put(0.5,3.5){1}
\put(0.5,6.3){2}
\put(0.5,7){y}
\color{red}
\put(2,6.5){$\G_{3}$}
\color{green}
\put(0.5,5){$\G_{4}$}
\color{blue}
\put(3.6,5){$\G_{2}$}
\put(2,3.4){$\G_{1}$}
\color{black}
\put(4.8,1.5){x}
\put(3.4,1.5){1}
\put(2,5){$\O$}
\put(1.2,6.4){D}
\put(1.2,3.5){A}
\put(3.5,6.4){C}
\put(3.5,3.5){B}
\thinlines
\qbezier[10](3.4,2)(3.4,2.8)(3.4,3.7)
\put(2,0){figure\ 1.1}
\end{picture}
\subsection{Example 2 : $\O=]0,1[\t]1,2[$,\ $l(x,y)=\sqrt{x}$, $\u(x,y)=(x,-y)$.}
\H Here, we have the same domain $\O$ as previously, therefore the sets $\G^-$, $\G^0$ and $\G^+$ are the same, but we have now a function 
$l$ that not belongs to $H^1(\O)$. As previously, from (\ref{2sz}) and $a=2$, in view of the boundary condition $z_{|\G^-}=0$, we derive $C=0$ and 
$$\forall (x,y)\in [0,1]\t [1,2],\ z(x,y)=\F{\sqrt{x}}{6}(4-y\sqrt{2y}).$$
We can verify that $z\in L^2(\O)$, but $z\notin H^1(\O)$. Clearly, the reason why is that $l\notin H^1(\O)$ and the problem (\ref{pbh2}) has no $H^1$ solution.\\[0.2cm]
\H In the following three examples, the function $\u\,.\,\mb{n}$ vanishes on $\ov{\G^-}$. In Example 3, the two assumptions of (\ref{4ctheun}) are verified and the solution $z$ is $H^1$, as expected by Theorem \ref{5theu}. On the contrary, in Examples 4 and 5, one of the two assumptions of (\ref{4ctheun}) is not verified (the first in Example 4 and the second in Example 5) and the solution $z$ is not $H^1$, thus proving the necessity of the two hypotheses (\ref{4ctheun}) in Theorem \ref{5theu}. 
 \subsection{Example 3 : $\O=\mr{triangle}(A(-\F{1}{2},\F{1}{2}),B(\F{1}{2},\F{1}{2}),C(\F{1}{2},\F{3}{2}))$,\\  $l(x,y)=1$, $\u(x,y)=(x,-y)$.}
 \begin{minipage}[b]{11cm}
 \H In this example, the set $\G^-$ is the line $]A,C[$ and the function $\u\,.\,\mb{n}_{|\G^-}$ vanishes at the endpoint $A$.\\
 \H We can verify that $\G^0=\emptyset$ and $\G^+=]A,B]\cup]B,C[$.
 From (\ref{2sz}), $a=\F{1}{2}$ and $l=1$ we derive
 $$\forall (x,y)\in\O,\ z(x,y)=1-2y+y\,C(xy).$$
 Considering the boundary condition, we have\\ $(z_{|\G^-}=0)\Longleftrightarrow (\forall(x,y)\in\G^-,\ C(xy)=2-\di\F{1}{y})$.\\ Setting $X=xy$, we have the 
 following equivalence:\\ $\forall (x,y)\in \G^-,\  \left\lbrace\begin{array}{l}y=x+1\\ -0.5<x<0.5\end{array}\right.
 \Longleftrightarrow\left\lbrace\begin{array}{l}y=\di\F{1+\sqrt{1+4X}}{2}\\ -0.25<X<0.75\end{array}\right.$
 \end{minipage}
\hspace*{1cm}
\unitlength=1cm
\begin{picture}(6,6)
\put(1,1.5){\vector(0,1){4.5}}
\put(-0.3,2){\vector(1,0){4}}
\color{blue}
\linethickness{0.04cm}
\put(0,3){\line(1,0){2}}
\put(2,3){\line(0,1){2}}
\color{red}
\put(2,5){\line(-1,-1){2}}
\color{black}
\put(0.5,5.7){y}
\color{red}
\put(0.1,3.7){$\G_{3}$}
\color{blue}
\put(1.2,2.4){$\G_{1}$}
\put(2.1,4){$\G_{2}$}
\color{black}
\put(0.6,4){1}
\put(0.35,5){1.5}
\put( 1.8, 1.6){0.5}
\put(-0.4, 1.6){-0.5}
\put(3,1.6){1}
\thinlines
\put(0.95,4){\line(1,0){0.1}}
\put(3.8,1.6){x}
\put(1.5,3.5){$\O$}
\put(0.55,1.55){O}
\put(-0.5,2.6){A}
\put(2,5.2){C}
\put(2.15,2.6){B}
\put(3,1.95){\line(0,1){0.1}}
\thinlines
\qbezier[10](2,2)(2,2.5)(2,3)
\qbezier[10](1,5)(1.5,5)(2,5)
\qbezier[10](0,2)(0,2.5)(0,3)
\put(1.7,0.4){figure\ 1.2}
\end{picture}\\[0.2cm]
Finally, we obtain the unique solution 
$$\forall(x,y)\in \O,\ z(x,y)=1-\di\F{2y}{1+\sqrt{1+4xy}}.$$
Indeed, $z\in L^2(\O)$ but we must verify that $z\in H^1(\O)$. We have
$$z'_x(x,y)=\di\F{4y^2}{(1+\sqrt{1+4xy})^2)}\di\F{1}{\sqrt{1+4xy}}.$$
For computing $\int\int_{\O}(z'_x)^2\,dxdy$, 
we make the substitution $\left\lbrace\begin{array}{l}X=xy\\
y=y\end{array}\right.$, the jacobian of which is $\di\F{1}{y}$. We obtain
$$\int\int_{\O}(z'_x)^2\,dxdy\le\int_{0.5}^{ 1.5}16y^3\,dy(\int_{y^2-y}^{0.5y}\di\F{1}{1+4X}\,dX)=\int_{0.5}^{1.5}4y^3(\ln(1+2y)-2\ln(2y-1))\,dy.$$
This last integral converges because $\int_{0.5}^{1.5}\ln(2y-1)\,dy$ is convergent in the neighbourhood of $0.5$. We can compute $\int\int_{\O}(z'_y)^2\,dxdy$
in the same way. Thus, we obtain that the solution $z$ belongs to $H^1(\O)$ and, therefore, the problem (\ref{pbh2}) is well-posed.\\
\H Note that, the fact that the function $\u\,.\,\mb{n}_{|\G^-}$ vanishes at the point $A$ leads to a discontinuity for the partial derivatives of the solution $z$ in $A$. 
However, the function $\u\,.\,\mb{n}_{|\G^-}$ has a simple root in $A$ and, moreover, $\u\,.\,\tho_-(A)=-\F{\sqrt{2}}{2}<0$ in $A$, where $\tho_-(A)$ is the unit tangent vector oriented towards $\G^-$. Thus, the assumptions (\ref{4ctheun}) of Theorem \ref{5theu} are verified, which explains that the solution $z$ is still in $H^1$.\\
\H In the two following examples, we change the function $\u$. In the Example 4, the function $\u\,.\,\mb{n}$ vanishes at the end point of $\G^-$ with an order two and  the solution $z$ is not $H^1$, which is consistent with the Theorem \ref{5theu}, since the assumption (\ref{4ctheun}) is not verified. In the Example 5, 
the function $\u\,.\,\mb{n}$ vanishes at the end point $A$ of $\G^-$ with an order one (simple root), but the  assumption (\ref{4ctheun}) is no longer verified, 
since the function $\u\,.\,\tho_-(A)$ is positive in $A$, and again, but for another reason, the solution $z$ is not $H^1$.
\subsection{Example 4 : $\O=\mr{triangle}(A(-2,\F{1}{3}),B(0,\F{1}{3}),C(0,1))$,\ $l(x,y)=1$, $\u(x,y)=(2xy+1,-y^2)$.}
 \begin{minipage}[b]{9cm}
 \H In this example, we change the function $\u$. We shall see that $\G^-$ is the line $]A,C[$ and, as in the example 3, the function $\u\,.\,\mb{n}_{|\G^-}$ va-\linebreak nishes at the endpoint $A$. However, contrary to the \linebreak example 3, \underline{the solution $z$ does not belong to $H^1(\O)$}. The reason why is that, contrary to the previous\linebreak example, the root in $A$ is a double root, as we will show below. We can verify that $\G^0=\emptyset$ and $\G^+=]A,B]\cup]B,C[$. 
     \end{minipage}
\hspace*{1cm}
\unitlength=1cm
\begin{picture}(6,5)
\put(4.15,1.5){\vector(0,1){3.5}}
\put(-0.5,2){\vector(1,0){5}}
\linethickness{0.04cm}
\color{red}
\qbezier(0,2.667)(2,3.333)(4,4)
\color{blue}
\qbezier(0,2.667)(2,2.667)(4,2.667)
\qbezier(4,2.667)(4,3)(4,4)
\color{black}
\put(3.5,5){y}
\put(4.5,1.7){x}
\put(3.6,1.6){O}
\put(0,3){A}
\put(3.6,4.2){C}
\put(3.6,2.9){B}
\put(-0.2,1.5){-2}
\put(4.2,2.7){$\F{1}{3}$}
\put(4.2,4){1}
\thinlines
\qbezier[10](0,2.667)(0,2.33)(0,2)
\put(1.7,0.4){figure\ 1.3}
\color{red}
\put(1.5,3.5){$\G^-$}
\color{blue}
\put(2.5,2.2){$\G^+$}
\end{picture}\\[0.2cm] 
We consider the following transport problem : find $z\in H^1(\O)$ satisfying
 \begin{equation}\label{pbh6}\left\lbrace
\begin{array}{l} z+\u\,.\,\n z=l\H \mr{in}\ \O,\\
z_{|\G^-}=0\end{array}\right. \end{equation}
\H Since, for all $(x,y)\in [A,C]$, $\u\,.\,\mb{n}=-\F{(x+2)^2}{\sqrt{10}},$ we obtain $\G^-=]A,C[$. As we saw previously, the solution z is expressed by
 (\ref{2sz2}). Next, in view of $X=xy^2+y$, we have the following equivalence$$  \left\lbrace\begin{array}{l}y=\F{x}{3}+1\\ -2\le x\le 0\end{array}\right.
 \Longleftrightarrow\left\lbrace\begin{array}{l}y=(\F{X-\F{1}{9}}{3})^{\F{1}{3}}+\F{1}{3}\\ \F{1}{9}\le X\le 1 \end{array}\right.$$
Then we derive $$\forall(x,y)\in \G^-,\ z(x,y)=0\Longleftrightarrow C(X)=e^{\F{1}{3}}-e^{\F{1}{y}}=e^{\F{1}{3}}-e^{\F{1}{(\F{X-\F{1}{9}}{3})^{\F{1}{3}}+\F{1}{3}}},$$
which allows us to compute the unique solution $z$ of Problem (\ref{pbh6}) :
\begin{equation}\label{3z6}
\forall (x,y)\in\O,\ z(x,y)=1-e^{\F{1}{\a(x,y)}-\F{1}{y}},
\end{equation}
with the function $\a$ defined in $\O$ by 
$$\forall (x,y)\in\O,\ \a(x,y)=(\F{xy^2+y-\F{1}{9}}{3})^{\F{1}{3}}+\F{1}{3}.$$
Since the domain $\O$ is below the segment$[AC]$ and since the branch of hyperbola $\left\lbrace\begin{array}{l}y=-\F{1}{x}\\ x<0\end{array}\right.$ is above the segment $[AC]$, we derive that, 
for all $(x,y)\in\ov{\O},\ 0<xy+1\le 1$. The function $(x,y)\mapsto xy+1$ is continuous on the compact $\ov{\O}$, therefore there exists $m_0>0$ such that 
$\forall (x,y)\in \O,\ xy+1\ge m_0$. Note that $m_0\le (-2)\F{1}{3}+1=\F{1}{3}$, which gives 
$$\forall (x,y)\in \O,\ m_0\le xy+1\le 1,\ \mr{with}\ 0<m_0\le\F{1}{3}.$$ Hence, we obtain
\begin{equation}\label{3ma6}
\forall (x,y)\in \O,\ \F{m_0}{3}\le\F{1}{3}(1-(1-3m_0)^{\F{1}{3}}))\le \a(x,y)\le 1.\end{equation}
Let us show that $z'_x$ does not belong to $L^2(\O)$. Considering (\ref{3z6}), we compute
$$\forall(x,y)\in\O,\ z'_x(x,y)=\F{\a'_x(x,y)}{(\a(x,y))^2}e^{\F{1}{\a(x,y)}-\F{1}{y}},$$
with $$\a'_x(x,y)=\F{y^2}{9(\F{xy^2+y-\F{1}{9}}{3})^{\F{2}{3}}}.$$
From (\ref{3ma6}), we derive
$$\forall (x,y)\in\O,\ |z'_x(x,y)|\ge\F{1}{3^4\,e^2}\F{1}{(\F{xy^2+y-\F{1}{9}}{3})^{\F{2}{3}}}. $$
Using this estimation yields
$$\int\int_{\O}(z'_x(x,y))^2\,dxdy\ge \F{1}{3^8\,e^4}\int_{\F{1}{3}}^1 dy\left(\int_{3(y-1)}^0\F{1}{(\F{xy^2+y-\F{1}{9}}{3})^{\F{4}{3}}}\,dx\right).$$
Making the substitution $\left\lbrace\begin{array}{l}X=xy^2+y\\
Y=\F{1}{y}\end{array}\right.$, the jacobian of which is -1, we obtain
$$\int\int_{\O}(z'_x(x,y))^2\,dxdy\ge \F{1}{3^8\,e^4}\int_1^3\,dY
\left(\int_{3(\F{1}{Y}-\F{1}{3})^3+\F{1}{9}}^{\F{1}{Y}}\F{1}{\left(\F{X-\F{1}{9}}{3}\right)^{\F{4}{3}}}\,dX\right)
\ge\F{1}{3^5\,e^4}\int_1^3\F{1}{3-Y}\,dY-\F{2^{\F{2}{3}}}{3^5\,e^4}.$$
$$\mr{Since}\ \int_1^3\F{1}{3-Y}\,dY=+\iy,\ \mr{we\ obtain}\ \int\int_{\O}(z'_x(x,y))^2\,dxdy=+\iy.\hspace*{3.9cm}$$
Finally, the solution $z$ of the example 4, contrary to the previous example, does not belong to $H^1(\O)$ and, therefore, the problem (\ref{pbh6}) is not well-posed.\\
\subsection{Example 5 : $\O=\mr{triangle}(A(-\F{4}{3},\F{2}{3}),B(-\F{11}{6},\F{1}{6}),C(-\F{4}{3},\F{1}{2}))$,\\ $l(x,y)=1$, $\u(x,y)=(-2xy-1,y^2)$.}
\hspace*{2cm}\unitlength=2cm
\begin{picture}(6,4)
\put(4.15,1.5){\vector(0,1){2.8}}
\put(-0.5,2){\vector(1,0){5}}
\linethickness{0.02cm}
\color{red}
\qbezier(1.333,3.333)(0.833,2.833)(0.667,2.667)
\qbezier(0.333,2.3333)(0.833,2.6667)(1.3333,3)
\color{blue}
\qbezier(0.667,2.667)(0.517,2.517)(0.333,2.333)
\qbezier(1.333,3)(1.333,3.2)(1.333,3.333)
\color{black}
\put(3.9,4.2){y}
\put(4.5,1.7){x}
\put(3.7,1.7){O}
\put(0.15,1.7){$-\F{11}{6}$}
\put(0.55,1.7){$-\F{5}{3}$}
\put(1.15,1.7){$-\F{4}{3}$}
\put(4.2,4){1}
\put(1.9,1.8){-1}
\put(1.2,3.5){A}
\put(0,2.2){B}
\put(1.4,2.8){C}
\put(0.45,2.7){D}
\put(4.1,2.3){$\F{1}{6}$}
\put(4.1,2.625){$\F{1}{3}$}
\put(4.1,2.95){$\F{1}{2}$}
\put(4.1,3.3){$\F{2}{3}$}
\thinlines
\qbezier[7](0.333,2)(0.333,2.16)(0.333,2.333)
\qbezier[10](0.667,2)(0.667,2.333)(0.667,2.667)
\qbezier[20](1.333,2)(1.333,2.5)(1.333,3)
\qbezier[40](4,2.333)(2,2.333)(0.333,2.333)
\qbezier[35](4,3)(2.6,3)(1.333,3)
\qbezier[35](4,3.333)(2.6,3.333)(1.333,3.333)
\qbezier[35](4,2.667)(2.3,2.667)(0.667,2.667)
\qbezier(4.07,4)(4.09,4)(4.11,4)
\qbezier(2,1.98)(2,2)(2,2.02)
\put(1.7,0.8){figure\ 1.4}
\color{red}
\put(0.6,3){$\G_1^-$}
\put(0.9,2.45){$\G_2^-$}
\color{blue}
\put(0.2,2.5){$\G_1^+$}
\put(1.4,3.1){$\G_2^+$}
\end{picture}\\[0.2cm]
\vspace*{-14cm}\\
\hspace*{2cm}\unitlength=8cm
\begin{picture}(3,3)
\linethickness{0.02cm}
\color{green}
\qbezier(0.64334,0.54)(0.6 ,0.6)(0.667,0.667)
\qbezier(0.667,0.667)(0.75,0.75)(1.333,1.18)
\put(1.15,0.98){$\gamma_1$}
\put(0.65,0.6){$\gamma_2$}
\color{red}
\qbezier(1.333,1.333)(0.833,0.833)(0.667,0.667)
\qbezier(0.333,0.3333)(0.833,0.6667)(1.3333,1)
\color{blue}
\qbezier(0.667,0.667)(0.517,0.517)(0.333,0.333)
\qbezier(1.333,1)(1.333,1.2)(1.333,1.333)
\color{black}
\put(1.32,1.36){A}
\put(0.28,0.26){B}
\put(1.35,0.95){C}
\put(0.6,0.65){D}
\put(1.15,1.1){$\O_1$}
\put(0.567,0.54){$\O_2$}
\put(0.8,0.7){$\O_3$}
\put(0.8,0.4){figure\ 1.5}
\color{red}
\put(0.85,1){$\G_1^-$}
\put(0.9,0.6){$\G_2^-$}
\color{blue}
\put(0.4,0.6){$\G_1^+$}
\put(1.4,1.1){$\G_2^+$}
\end{picture}
\H We can verify: $(x,y)\in (AB)\Leftrightarrow\ x-y+2=0$ and $(x,y)\in (BC)\Leftrightarrow\ 2x-3y+\F{25}{6}=0$.
 The set $\G^-$ is composed of two parts :
$\G^-_1\ \left \lbrace\begin{array}{l} x-y+2=0\\
 \F{1}{3}<y<\F{2}{3}\end{array}\right.$
 and $\G^-_2\ \left \lbrace\begin{array}{l} 2x-3y+\F{25}{6}=0\\
 \F{1}{6}<y<\F{1}{2}\end{array}\right..$
  The set $\G^+$ is composed of two parts :
$\G^+_1\ \left \lbrace\begin{array}{l} x-y+2=0\\
 \F{1}{6}<y<\F{1}{3}\end{array}\right.$
 and $\G^+_2\ \left \lbrace\begin{array}{l} x=-\F{4}{3}\\
 \F{1}{2}<y<\F{2}{3}\end{array}\right..$\\
 \H For all $\m=(x,y)\in \G_1^-$, $\u\,.\,\mb{n}(\m)=(3x+5)(x+1)$ and for all $\m=(x,y)\in \G_2^-$, $\u\,.\,\mb{n}(\m)=-9y^2+\F{25}{3}y-2$. Therefore, the function $\u\,.\,\mb{n}_{|\G}$ vanishes at the unique point $D(-\F{5}{3},\F{1}{3})$, with an order one with respect to the parameter of the line $(AB)$. Thus, we have 
 \begin{equation}\label{2ex5}
 \u\,.\,\mb{n}_{|\G}(D)=0,\ \F{\P\u}{\P\tho_-}\,.\,\mb{n}(D)\not=0\ \mr{and}\ (\u\,.\,\tho_-)_{|\G}(D)>0,\end{equation} where $\tho_-(\F{\sqrt{2}}{2},\F{\sqrt{2}}{2})$ is unit tangent vector, oriented towards $\G^-$.\\
 \H Setting $\left\lbrace\begin{array}{l}X=-xy^2-y\\
 Y=\F{1}{y}\end{array}\right.$, by technics analogous to the previous examples, we obtain the solutions of the equation $z+\u\,.\,\n z=l$
 $$Z(X,Y)=1+e^Y(C(X)-1)\Leftrightarrow z(x,y)=1+e^{\F{1}{y}}(C(-xy^2-y)-1),$$where $C$ is a function to be determined by the boundary conditions.\\
 \H Setting $\a(y)=X(y-2,y)=-y^3+2y^2-y$, we can verify
 $$z\,.\,\mb{n}_{|\G_1^-}=0\Longleftrightarrow\forall y\in ]\F{1}{3},\F{2}{3}[,\ C(\a(y))=1-e^{\F{1}{y}}\Longleftrightarrow\forall X\in ]-\F{4}{27},-\F{2}{27}[,\ C(X)=1- e^{\F{1}{\a^{-1}(X)}}.$$
 In the same way, setting $\beta(y)=X(\F{3}{2}y-\F{25}{12},y)=-\F{3}{2} y^3+\F{25}{12} y^2-y$, we can verify
 $$z\,.\,\mb{n}_{|\G_2^-}=0\Longleftrightarrow\forall y\in ]\F{1}{6},\F{1}{2}[,\ C(\beta(y))=1-e^{\F{1}{y}}\Longleftrightarrow\forall X\in ]-\F{1}{6},-\F{25}{216}[,\ C(X)=1- e^{\F{1}{\beta^{-1}(X)}}.$$
 \H Taking into account these boundary conditions, setting $y_1=\beta^{-1}(-\F{4}{27})$ and using a function $\a_1$, which is a restriction of the function $\a$, and functions $\a_2$ and $\a_3$, which are restrictions of the function $\beta$, we express the solution $z$ by splitting the domain $\O$ into three sub domains $\O_i$, $i=1,2,3$ :
 \begin{equation}\label{ex5sz}
 z_{|\O_i}=1-e^{\F{1}{y}-\di\F{1}{\a_i^{-1}(-xy^2-y)}},\end{equation}
 where $\O_1$ is defined by
 $$\O_1=\{(x,y)\in\O,\ y>\F{1}{3},\ -xy^2-y> -\F{4}{27}\}\ \mr{with}\ \begin{array}{ccc}
 &&\\
\a_1 :&[\F{1}{3},\F{2}{3}]\to &[-\F{4}{27},-\F{2}{27}]\\
 &\hspace*{0.5cm}y\hspace*{0.5cm}\longmapsto&-y^3+2y^2-y\end{array},$$
 where $\O_2$ is defined by
 $$\O_2=\{(x,y)\in\O,\ y<\F{1}{3},\ -xy^2-y> -\F{4}{27}\}\ \mr{with}\ \begin{array}{ccc}
 &&\\
\a_2 :&[\F{1}{6},y_1]\to &[-\F{4}{27},-\F{25}{216}]\\
 &\hspace*{0.5cm}y\hspace*{0.5cm}\longmapsto&-\F{3}{2}y^3+\F{25}{12}y^2-y\end{array}$$
 and where $\O_3$ is defined by
 $$\O_3=\{(x,y)\in\O,\ -xy^2-y< -\F{4}{27}\}\ \mr{with}\ \begin{array}{ccc}
 &&\\
\a_3 :&[y_1,\F{1}{2}]\to &[-\F{1}{6},-\F{4}{27}]\\
 &\hspace*{0.5cm}y\hspace*{0.5cm}\longmapsto&-\F{3}{2}y^3+\F{25}{12}y^2-y\end{array}.$$
 \H Note that the domains $\O_1$ and $\O_3$ are adjacent and are separated by the curve 
 $\gamma_1$ and the domains $\O_2$ and $\O_3$ are adjacent and are separated by the curve $\gamma_2$ (see the figure 1.5), where $\gamma_1$ and $\gamma_2$ are defined by
   $$\gamma_1\ \left\lbrace\begin{array}{l}-x y^2-y=-\F{4}{27}\\
 \F{1}{3}<y<\F{3+\sqrt{3}}{8}\end{array}\right.,\ \gamma_2\ \left\lbrace\begin{array}{l}-x y^2-y=-\F{4}{27}\\
 y_1<y<\F{1}{3}\end{array}\right..$$
 \H Considering the expressions of the solution given by (\ref{ex5sz}), we obtain 
 $$\forall y\in ]\F{1}{3},\F{3+\sqrt{3}}{8}[,\ (z_{|\O_1}-z_{|\O_3})_{|\gamma_1}(y)=(e^{-\F{1}{y_1}}-e^{-3})e^{\F{1}{y}},$$
 which implies that the solution $z$ is discontinuous on the curve $\gamma_1$. Computing the gradient of the solution $z$ yields
 $$\n z=\widetilde{\n z_{|\O_1}}+\widetilde{\n z_{|\O_2}}+\widetilde{\n z_{|\O_3}}-\delta_{\gamma_1},$$
 where the wide tildes denote the extensions by zero and where the distribution $\delta_{\gamma_1}$ is defined by
 $$\forall\vb\in (\mathcal{D}(\O))^2,\ <\delta_{\gamma_1},\vb>=\int_{\gamma_1}(z_{|\O_1}-z_{|\O_3})\vb\,.\,\mb{n_1}\,ds,$$
 where $\mb{n_1}$ is the unit exterior normal vector to the boundary of the domain $\O_1$. Finally, since the distribution $\delta_{\gamma_1}$ does not belongs to $L^2(\O)$, we obtain that the solution $z$ of the example 5 does not belong to $H^1(\O)$ and, therefore, the problem (\ref{pbh6}) is not well-posed.\\
 \H For explaining further in details, in view of (\ref{2ex5}), the function $\u\,.\,\mb{n}$ vanishes at the boundary point $D$ of $\Gamma^-_1$ with an order one with respect to the parameter of the line $(AB)$, but,  since we have $(\u\,.\,\tho_-)_{|\G}(D)>0$, the solution $z$ in the neighborhood of $D$ on the $\O_2$ side depends of the boundary condition on $\Gamma_2^-$, which is far from $D$. This means that we cannot localize the transport problem in a neighborhood of the boundary point $D$ and, therefore, we cannot apply the technics of the proof of Theorem \ref{5theu}. Thus, the assumption $\u(\m)\,.\,\tho_-(\m)<0$ of Theorem \ref{5theu} is not only a technical assumption, but a basic assumption as well as the other assumption $\F{\P\u}{\P\tho_-}(\m)\,.\,\mb{n_-}(\m)\not= 0$ of (\ref{4ctheun}).
     \section{Transport equations in $H^1$ when $\u\,.\,\mb{n}$ does not \\ vanish on $\ov{\G^-}$}
\H Let us recall the following problem studied in [3]. Let $\O$ be a bounded domain of $\R^2$  and $\G^-$ be defined by (\ref{1dg-}), 
verifying (\ref{1dg0g1}): for $\u$ in $H^1(\O)^d$, with $\d\u=0$, $l$ in $L^2(\O)$ and $\W$ in 
$\R^*$, find $z$ in $L^2(\O)$ such that
\begin{equation}\left\lbrace\begin{array}{ll}
z+\W\,\u\,.\,\n z=l\H &\mr{in}\H \O\\
(z\u)\,.\,\mb{n}=0 &\mr{on}\H \G^-.\end{array}\right.\label{4eqt1}
\end{equation}
\H The main result is given by Theorem 3.3 in [3], which gives the existence and the uniqueness of solution in $L^2(\O)$ in the case 
where $\O$ is a Lipschitz-continuous domain of $\R^d$. Now, we are interested by $H^1$ solutions in the two dimensions case. In order to find 
 $H^1$ solutions, we assume that $\O$ is a bounded polygon, we suppose that $\u$ belongs to $W^{1,\iy}(\O)^2$ and we shall impose another boundary condition.\\ 
  \H Thus, we are led to study the following problem: let $\O$ be a bounded polygon,  
 for $\u$ in  $U\cap W^{1,\iy}(\O)^2$, where $U$ is defined by (\ref{dU}), $l$ in $H^1(\O)$ and $\W$ in $\R^*$, 
 find $z$ in $H^1(\O)$ such that
 \begin{equation}\left\lbrace\begin{array}{ll}
z+\W\,\u\,.\,\n z=l\H &\mr{in}\H \O\\
z=0 &\mr{on}\H \G^-.\end{array}\right.\label{4eqt2}
\end{equation}
  \H Let $\O$ be a bounded polygon. We begin to establish a result of existence and uniqueness in the particular where $l$ vanishes on $\G^-$.
 \begin{theo}\label{4theu1} Let $\O$ be a bounded polygon, $\G^-$ be defined by (\ref{1dg-}), verifying 
(\ref{1dg0g1}) and $U$ be defined by (\ref{dU}). For all $\u$ in 
 $U\cap W^{1,\iy}(\O)^2$ such that 
 \begin{equation}\label{4ctheub}\p\n\u\p_{L^{\iy}(\O)}\le \di\F{1}{2|\W|},\end{equation} 
 all $l$ in $H^1(\O)$ such that $l_{|\G^-}=0$
  and all real number $\W$ in $\R^*$, the transport problem (\ref{4eqt2}) has a unique solution $z$ in $H^1(\O)$.
\end{theo}
\tb{Proof.} Formally, $\n z$ satisfies
$$\n z+\W\,\u\,.\,\n(\n z)=\n l-\W\,\n\u\,.\,\n z.$$
\H Let us define a sequence $(\mb{F}_n)$ of functions $\mb{F}_n\in X_{\u}(\G^-)^2$, $n\in \N$, by recurrence, where 
$X_{\u}(\G^-)$ is defined by (\ref{1dXug0}). We set $\mb{F}_0=\mb{0}$ and assume that the function $\mb{F}_n\in X_{\u}(\G^-)^2$ is given for $n\in\N$. Then, applying Theorem \ref{2teeth}, we define each component $F_{n+1,1}$ and $F_{n+1,2}$ of $\mb{F}_{n+1}$ as the unique solution of a transport equation from the type (\ref{4eqt1}), of such so that we define $\mb{F}_{n+1}\in X_{\u}(\G^-)^2$ as the unique solution of the transport equation
\begin{equation}\left\lbrace\begin{array}{ll}
\mb{F}_{n+1}+\W\,\u\,.\,\n \mb{F}_{n+1}=\n l-\W\,\n\u\,.\,\mb{F}_{n}\H &\mr{in}\H \O\\
(\mb{F}_{n+1}\,\u)\,.\,\mb{n}=0 &\mr{on}\H \G^-.\end{array}\right.\label{4eqt3}
\end{equation}
Since $\mb{F}_{n+1}$ belongs to $X_{\u}(\G^-)^2$, the basic result of Proposition (\ref{2pgreenfg-}) implies
$$\int_{\O}(\W\,\u\,.\,\n F_{n+1,i})F_{n+1,i}\,d\mb{x}\ge 0,\ \mr{for}\ i=1,2.$$
Then, taking the scalar product of both sides of the first equation of (\ref{4eqt3}) with $\mb{F}_{n+1}$ yields
$$\p\mb{F}_{n+1}\p_{L^2(\O)}^2\le (\n l,\mb{F}_{n+1})-\W\,(\n\u\,.\,\mb{F}_{n},\mb{F}_{n+1}).$$
Hence, we derive
$$\p\mb{F}_{n+1}\p_{L^2(\O)}\le \p\n l\p_{L^2(\O)}+|\W|\p\n\u\p_{L^{\iy}(\O)}\p\mb{F}_{n}\p_{L^2(\O)}.$$
In view of the bound (\ref{4ctheub}), we obtain
$$\p\mb{F}_{n+1}\p_{L^2(\O)}\le \p\n l\p_{L^2(\O)}+\F{1}{2}\p\mb{F}_{n}\p_{L^2(\O)},$$
which implies, by a recurrence argument, that $\mb{F}_{n}$ is uniformly bounded in $L^2(\O)$ and 
$\forall n\in \N$,
\begin{equation}\label{4bFn}
\p\mb{F}_{n}\p_{L^2(\O)}\le 2\p\n l\p_{L^2(\O)}.\end{equation}
Owing to (\ref{4bFn}), $\u\,.\,\n \mb{F}_{n+1}$ is also uniformly bounded in $L^2(\O)$. Therefore we can pass to the limit in the first equation of 
(\ref{4eqt3}) and there exists a function $\mb{F}\in L^2(\O)^2$ such that
\begin{equation}\label{4eqtd}\mb{F}+\W\,(\u\,.\,\n\mb{F}+\n\u\,.\,\mb{F})=\n l.\end{equation}
Let us set $z=l-\W\u\,.\mb{F}$. From the previous equation, we derive $\mb{F}=\n z$ and we obtain $z=l-\W\u\,.\,\n z$, which gives that $z$ is solution of the first equation of (\ref{4eqt2}).\\
\H Next, from Green's formula (\ref{1green2}) and $(\mb{F}_{n+1}\u)\,.\,\mb{n}_{|\G^-}=\mb{0},$ we derive
 $\forall\vfy\in W^{1,r}(\O)^2$, with $\vfy_{|\G^{0,+}}=\mb{0}$,  
 $$(\mb{F}_{n+1}\u,\n\vfy)+(\vfy\u,\n \mb{F}_{n+1})= <(\mb{F}_{n+1}\u)\,.\,\mb{n},\vfy>_{\G^-}=0.$$
 Using the above convergence, we can pass to the limit and we obtain
 $$\forall\vfy\in W^{1,r}(\O)^2,\ \mr{with}\ \vfy_{|\G^{0,+}}=0,\ (\mb{F}\,\u,\n\vfy)+(\vfy\u,\n \mb{F})=0,$$
 which implies, with again the Green's formula (\ref{1green2}), $<(\mb{F}\,\u)\,.\,\mb{n},\vfy>_{\G^-}=0$.
 Thus, we obtain $(\mb{F}\,\u)\,.\,\mb{n}_{|\G^-}=0$, that is to say,
 \begin{equation}\label{4cbd}(\n z\,\u)\,.\,\mb{n}_{|\G^-}=0.\end{equation}
 Hence, we can use a density result of [3](Corollary 2.11, page 1012): since, for $i=1,2$, $\di\F{\P z}{\P x_i}$ belongs to $X_{\u}(\G^-)$, there exist two sequences $(\vfy_{1,n})$ et $(\vfy_{2,n})$ such that, for $i=1,2$, $\vfy_{i,n}\in \mathcal{D}(\ov{\O},\G^-)$ and 
 $$\lim_{n\to +\iy}\vfy_{i,n}=\F{\P z}{\P x_i}\H \mr{strongly\ in}\H X_{\u}(\G^-), $$ 
 where $X_{\u}(\G^-)$ is defined in (\ref{1dXug0}). Setting $\vb_n=\left(\begin{array}{c}\vfy_{1,n}\\ \vfy_{2,n}\end{array}\right)$, 
 from the above convergence and the regularity of $\u$ we derive
 $$\lim _{n\to +\iy}\u\,.\,\vb_n=\u\,.\,\n z\H\mr{strongly\ in}\H L^2(\O).$$ 
 Noting that $$\n(\u\,.\,\n z)=\left (\begin{array}{c}\di\F{\P\u}{\P x_1}\,.\,\n z+\u\,.\,\n(\F{\P z}{\P x_1})\\
 \di\F{\P\u}{\P x_2}\,.\,\n z+\u\,.\,\n(\F{\P z}{\P x_2})\end{array}\right),$$ the convergences in $X_{\u}(\G^-)$ give, for $i=1,2$,
 $$\lim _{n\to +\iy}(\u\,.\,\n \vfy_{i,n})=\u\,.\,\n(\F{\P z}{\P x_i})\H\mr{strongly\ in}\H L^2(\O),$$
 $$\lim _{n\to +\iy}(\F{\P\u}{\P x_i}\,.\,\vb_n)=\F{\P\u}{\P x_i}\,.\,\n z\H\mr{strongly\ in}\H L^2(\O).$$
 These convergences imply
 $$\lim _{n\to +\iy}\n(\u\,.\,\vb_n)=\n(\u\,.\,\n z).$$
 Thus, we obtain that $$\lim_{n\to +\iy}(\u\,.\,\vb_n)=\u\,.\,\n z \H\mr{strongly\ in}\H H^1(\O).$$
  In view of $\vb_{n|\G^-}=\mb{0}$, we obtain
  $$(\u\,.\,\n z)_{|\G^-}=0.$$
 Considering that $z+\W\,\u\,.\,\n z=l$ and $l_{|\G^-}=0$, we obtain 
  $$z_{|\G^-}=0.$$
  Thus, we have proven the existence of solution 
 for the transport problem (\ref{4eqt2}).\\
 \H Concerning the uniqueness, let us consider $z\in H^1(\O)$ solution of the problem 
 \begin{equation}\left\lbrace\begin{array}{ll}
z+\W\,\u\,.\,\n z=0\H &\mr{in}\H \O\\
z=0 &\mr{on}\H \G^-.\end{array}\right.\label{4eqt4}
\end{equation}
For proving the uniqueness of solution of Problem (\ref{4eqt2}), we must show that necessarily $z=0$.
Taking the scalar product in $L^2(\O)$ of the previous equation by $z$ yields
$$\p z\p_{L^2(\O)}+\W(\u\,.\n z, z)=0.$$
Since $z$ belongs to $X_{\u}(\G^-)$, Proposition (\ref{2pgreenfg-}) implies $\W(\u\,.\n z, z)\ge 0$ and we derive
$$\p z\p_{L^2(\O)}\le 0.$$
This gives $ z=0$, which gives the uniqueness of solution of Problem (\ref{4eqt2}).
\hfill$\diamondsuit$.\\[0.3cm]   
 \H Now, we do not assume that $l$ vanishes on $\G^-$. If $\m$ belongs to $\ov{\G^{-}}$ and does not belong to $\ov{\G^{+,0}}$, we denote by
 \begin{equation}\label{4dnm1}
 \mb{n}_{-}(\m)\ \mr{the\ unit\ exterior\ normal\ vector\ to}\ \ov{\G^{-}}\ \mr{in}\ \m\end{equation}
(one or other of the two unit exterior normal vectors if $\m$ is a vertex of the polygon).  If $\m$ belongs to $\ov{\G^{-}}\cap\ov{\G^{+,0}}$, then $\m$ is the common endpoint of two adjacent straight segments $\gamma_+$ and $\gamma_-$ such that $\gamma_+\subset \ov{\G^{+,0}}$ and $\gamma_-=[\m,\mb{m_-}]\subset \ov{\G^-}$ with $\m\not=\mb{m_-}$. We denote by
 \begin{equation}\label{4dnm2}
 \mb{n}_{-}(\m)\ \mr{the\ unit\ exterior\ normal\ vector\ to}\ \gamma_-,\end{equation}  and by 
 \begin{equation}\label{4dtom}
 \tho_{-}(\m)\ \mr{the\ unit\ tangent\ vector}\ \F{1}{\p \m\mb{m_-}\p}\m\mb{m_-}.\end{equation}  
 \H First, we assume that the normal component of the velocity does not vanish on $\ov{\Gamma^-}$. Since $\u\,.\,\mb{n}$ is continuous on the sides of the polygon 
 $\O$, this implies that the end points of $\Gamma^-$ are vertices of the polygon. The following theorem gives assumptions implying existence and uniqueness for problem (\ref{4eqt2}). 
\begin{theo}\label{4theu2} Let $\O$ be a bounded polygon, $\G^-$ be defined by (\ref{1dg-}), verifying 
(\ref{1dg0g1}) and $U$ be defined by (\ref{dU}). For all $\u$ in 
 $U\cap W^{1,\iy}(\O)^2$ such that 
 \begin{equation}\label{4ctheu}\p\n\u\p_{L^{\iy}(\O)}\le \di\F{1}{2|\W|}\end{equation} and such that
 \begin{equation}\label{4cun}
 \forall \mb{m}\in\ov{\Gamma^-},\ \u(\m)\,.\,\mb{n}_-(\m)\not=0,\end{equation}
 where $\mb{n}_-(\m)$ is defined by (\ref{4dnm1}) or  by (\ref{4dnm2}), all $l$ in $H^1(\O)$
  and all real number $\W$ in $\R^*$, the transport problem (\ref{4eqt2}) has a unique solution $z$ in $H^1(\O)$.
\end{theo}
\tb{Proof.} Since the end points of $\G^-$ are vertices, we have 
\begin{equation}\label{4dGj-}
\ov{\G^-}=\bigcup_{j\in J}\G_j,\end{equation}
where the sets $\G_j$ are sides of the polygon $\O$. Since $\u$ is continuous on $\P\O$, for all $j\in J$, we denote $\eta_j=\min\limits_{m\in \G_j}(|\u(\m)\,.\,\mb{n}_j|)$. From (\ref{4cun}), we derive that, 
for all $j\in J$, $\eta_j>0$, which implies that $(\di\F{l}{\u\,.\,\mb{n_j}})_{|\G_j}$ belongs to $H^{\F{1}{2}}(\G_j)$. 
So, there exists $z_0$ in $H^2(\O)$ verifying, for all $j\in J$, 
 \begin{equation}\left\lbrace\begin{array}{l}
(\di\F{\P z_0}{\P n})_{|\G_j}=(\F{l}{\W\,\u\,.\,\mb{n_j}})_{|\G_j}\\
z_{0|\G_j}= 0.\end{array}\right.\label{4lghlp}\end{equation}
Hence, we derive that
\begin{equation}\label{4z0l}
(z_0+\W\,\u\,.\,\n z_0)_{|\G^-}=l_{|\G^-}.\end{equation}
Next, applying Theorem \ref{4theu1}, let $z^*\in H^1(\O)$ be the unique solution of the problem
$$\left\lbrace\begin{array}{ll}
z^*+\W\,\u\,.\,\n z^*=l-z_0-\W\,\u\,.\,\n z_0\H &\mr{in}\H \O\\
z^*=0 &\mr{on}\H \G^-\end{array}\right..$$
Then, $z=z^*+z_0$ verifies $z+\W\,\u\,.\,\n z=l$ and $z_{|\G^-}=0$. Thus, we have proven the existence of solution 
 for the transport problem (\ref{4eqt2}). We prove the uniqueness in the same way as in the previous theorem.\hfill$\diamondsuit$.\\[0.2cm]
\section{Transport equations in $H^1$ when $\u\,.\,\mb{n}$ vanishes on $\ov{\G^-}$}
\H We assume that $\O$ is a bounded convex polygon, but the fact that the normal component of the velocity can vanish on the boundary 
introduces a singularity at the end points of $\G^-$ and we will be forced to make assumptions at the end points of $\G^-$, as we could expect from the examples of the Section 2. We denote by
\begin{equation}\label{4dS}
S\ \mr{the\ set\ of\ the\ vertices\ of\ the\ polygon}\ \O\end{equation}
and let the set $E$ be defined by
\begin{equation}\label{4d-+0}
E=\{\m\in \ov{\Gamma^-}\cap\ov{\Gamma^{+,0}},\ \u(\m)\,.\,\mb{n}_-(\m)=0\},\end{equation}
where $\mb{n}_-(\m)$ is defined by (\ref{4dnm2}).
Note that, in view of the assumption (\ref{1dg0g1}), the set $E$ is finite. In addition, we make the assumption that the velocity $\u$ is such that
\begin{equation}\label{4hs}
\{\m\in\ov{\Gamma^-},\ \u(\m)\,.\,\mb{n}_-(\m)=0\}\subset E,\end{equation}
 which means that $\u\,.\,\mb{n}$ does not vanish in a point located in the interior of $\ov{\G^-}$.\\
\H The next theorem, which is the main result of the paper, gives assumptions implying existence and uniqueness for problem (\ref{4eqt2}), in the case where the normal component of the velocity vanishes on the boundary. Note that, the first assumption of (\ref{4ctheun}) means that the function $\u\,.\,\mb{n}$ must have only simple roots at the end points of $\G^-$, which seems consistent with the previously studied examples. At first glance, the second assumption of (\ref{4ctheun}) 
seems to be a technical assumption, related to the method used in the proof of Theoren \ref{5theu}. Indeed, we need this assumption, in the proof of the theorem, 
probably because, in \\ \pagebreak \\the case where $\u(\m)\,.\,\tho_-(\m)>0$, it does not seem possible to localize the problem around the points of the set $E$ : on either side of the point where $\u\,.\,\mb{n}$ vanishes, the expressions of the solution $z$ are determined by boundary conditions located in two different places of the boundary, which leads to a discontinuity of the solution $z$, see Example 5. In fact, as it appears in Example 5, this second assumption seems necessary to obtain a solution $z$ in $H^1$.
\begin{theo}\label{5theu} Let $\O$ be a bounded convex polygon, $\G^-$ be defined by (\ref{1dg-}), verifying 
(\ref{1dg0g1}) and  $U$ be defined by (\ref{dU}).  
For all $\u$ in 
 $U\cap W^{1,\iy}(\O)^2$, verifying (\ref{4hs}),  such that
 \begin{equation}\label{4ctheut}\p\n\u\p_{L^{\iy}(\O)}\le \di\F{1}{2|\W|}\end{equation}
 and such that 
 \begin{equation}\label{4ctheun}
 \forall \m\in E,\ \F{\P\u}{\P\tho_-}(\m)\,.\,\mb{n_-(\m)}\not= 0\ \mr{and}\ \u(\m)\,.\,\tho_-(\m)<0,\end{equation}
 where $\mb{n_-(\m)}$ (respectively $\tho_-(\m)$, $E$) is defined by (\ref{4dnm2}) (respectively (\ref{4dtom}), (\ref{4d-+0})), all $l$ in $H^1(\O)$
  and all real number $\W$ in $\R^*$, the transport equation (\ref{4eqt2}) has a unique solution $z$ in $H^1(\O)$.
\end{theo}
\tb{Proof.} Let us split up $\ov{\G^-}$ into straight segments as
\begin{eqnarray}
\ov{\G^-}=\bigcup\limits_{j=1}^q\gamma_j,\ \gamma_j\cap\gamma_k=\emptyset\ \mr{if}\ k\notin\{j-1,j,j+1\},\hspace*{1cm}\nonumber\\
\gamma_j\cap\gamma_k=\emptyset\ \mr{or}\ \gamma_j\cap\gamma_k\in S\ \mr{if}\ k\in\{j-1,j+1\},\ 1\le j\le q,\ 0\le k\le q+1\label{4decompg-}\end{eqnarray}
with the convention $\gamma_{q+l}=\gamma_l$ for $l=0,1$, and let $\mu_0>0$ be defined by
\begin{equation}\label{4dmu0}
\mu_0=\min_{\stackrel{\gamma_j\cap\gamma_k=\emptyset}{
1\le j,k\le q}} d(\gamma_j,\gamma_k),\end{equation}
where $d(.,.)$ is the euclidian distance in $\R^2$. Then, for $0<\mu\le \F{1}{2}\mu_0$, in order to localize around the sets $\gamma_j$, let us define the functions $(\theta_{j,\mu})_{1\le j\le q}\in\mathcal{D}(\R^2)$ by
\begin{equation}\label{4dthjmu}
\forall\x\in\R^2,\ \theta_{j,\mu}(\x)=\left\lbrace \begin{array}{l} 1\ \mr{if}\ d(\x,\gamma_j)\le\F{1}{2}\mu\\
0\ \mr{if}\ d(\x,\gamma_j)\ge \mu.\end{array}\right.\end{equation}and, $\forall \x\in\R^2$, $\theta_{q+1,\mu}(\x)=0$.         
Setting, for $1\le j\le q$ and $0<\mu\le\F{1}{2}\mu_0$, 
\begin{equation}\label{4deflmuj}
l_{j,\mu}=\theta_{j,\mu}(1-\theta_{j+1,\mu})l\end{equation} and 
\begin{equation}\label{4deflmu}
l_{\mu}=(1 -\sum_{j=1}^q\theta_{j,\mu}(1-\theta_{j+1,\mu}))l,\end{equation}
where $l$ is the right hand side of the transport equation, we obtain
\begin{equation}\label{4decpl}
l=l_{\mu}+\sum_{j=1}^q l_{j,\mu}\end{equation}
and we can verify that
\begin{equation}\label{4plmu}
\forall\x\in\G^-,\ l_{\mu}(\x)=0.\end{equation}
 From the development of $l$ given by (\ref{4decpl}), we derive $q+1$ problems, constructed from (\ref{4eqt2}) by substituting $l_{\mu}$, $l_{j,\mu}$, $1\le j\le q$, to $l$. First, the problem $(P_{\mu})$ : find $z$ in $H^1(\O)$ such that
 \begin{equation}(P_{\mu})\\ \left\lbrace\begin{array}{ll}
z+\W\,\u\,.\,\n z=l_{\mu}\H &\mr{in}\H \O\\
z=0 &\mr{on}\H \G^-\end{array}\right.\label{4eqtmu}\end{equation}
and, second, the problems $(P_{j,\mu})_{1\le j\le q}$ : 
find $z$ in $H^1(\O)$ such that
 \begin{equation}(P_{j,\mu})\ \left\lbrace\begin{array}{ll}
z+\W\,\u\,.\,\n z=l_{j,\mu}\H &\mr{in}\H \O\\
z=0 &\mr{on}\H \G^-.\end{array}\right.\label{4eqtjmu}
\end{equation}
Note that, because of the linearity, the solution of the problem (\ref{4eqt2}) will be the sum of the solution of the problem $(P_{\mu})$ and the solutions of problems $(P_{j,\mu})_{1\le j\le q}$.\\
\H In view of (\ref{4plmu}), applying Theorem \ref{4theu1}, we derive that 
\begin{equation}\label{4zmu}\mr{the\ problem}\ (P_{\mu})\ \mr{has\ a\ unique\ solution}\ z_{\mu}\in H^1(\O).
\end{equation}
\H Next, we have to solve the problems $(P_{j,\mu})_{1\le j\le q}$. We denote by $\mb{n_j}$ the exterior unit normal  vector of the side of the polygon which contains $\gamma_j$ and, for $i=-1,1$, by  $\mb{S^i_j}$ the end points of $\gamma_j$, with the convention that, if $\gamma_j\cap\gamma_{j+i}\not=\emptyset$ for $i=-1$ or $i=1$, then $\gamma_j\cap\gamma_{j+i}=\{\mb{S^i_j}\}$. Note that, for each point $\mb{S^i_j}$, $i=-1,1$, $1\le j\le q$, we have four possibilities : $\mb{S^i_j}\in \gamma_{j+i}$, $\mb{S^i_j}\notin\gamma_{j+i}$ with $\mb{S^i_j}\notin E$, $\mb{S^i_j}\notin\gamma_{j+i}$ with $\mb{S^i_j}\in (E\cap S)$, 
$\mb{S^i_j}\notin\gamma_{j+i}$ with $\mb{S^i_j}\in (E\cap S^c)$, where $S^c$ is the complementary set of $S$ in $\R^2$. We shall not consider all the cases, because there are similar cases, but we shall study some cases, which will be models for the other cases. Note that , for $i=1,2$, if $\mb{S^i_j}$ is not a vertex of the polygon, then $\u\,.\,\mb{n_j}(\mb{S^i_j})=0$, that is to say $\mb{S^i_j}\in E$. \\ 
1) \tb{First case}: $\mb{S^i_j}\in \gamma_{j+i}$, $i=-1,1$.\\
Note that, in view of (\ref{4decompg-}) and (\ref{4hs}), $\mb{S^i_j}\in S$ and $\u(\mb{S^i_j})\,.\,\mb{n_{j+i}}\not=0$, $\u(\mb{S^i_j})\,.\,\mb{n_j}\not=0$, for $i=-1,1$. Moreover, $l_{j,\mu}=0$ on $\gamma_k$, for $k\notin \{j-1,j\}$.
Since $\u(\mb{S^{-1}_j})\,.\,\mb{n_{j-1}}\not=0$, there exist a real number $\mu_1>0$ such that, for all $\x$ verifying $d(\mb{S^{-1}_j},\x)\le\mu_1$, we have
$\u(\mb{x})\,.\,\mb{n_{j-1}}\not=0$. Then, with the notation 
$$\gamma_{j-1,1}=\{\x\in\gamma_{j-1},\ d( \mb{S^{-1}_j},\x)\le \mu_1\},\ \gamma_{j-1,2}=\{\x\in\gamma_{j-1},\ d( \mb{S^{-1}_j},\x)> \mu_1\},$$ taking
\begin{equation}\label{4cmu01}0<\mu\le\min(\F{1}{2}\mu_0,\mu_1),\end{equation} in the same way as in the proof of Theorem \ref{4theu2}, there exists $z_{0,j,\mu}$ in $H^2(\O)$ verifying,  
$$\left\lbrace\begin{array}{l}
(\di\F{\P z_{0,j,\mu}}{\P n})_{|\gamma_j}=(\F{l_{j,\mu}}{\W\,\u\,.\,\mb{n_j}})_{|\gamma_j}\\
z_{0,j,\mu|\gamma_j}= 0\end{array},\right.\hspace*{2.15cm}$$
\begin{equation}\label{4z0jj-1}
\left\lbrace\begin{array}{l}
(\di\F{\P z_{0,j,\mu}}{\P n})_{|\gamma_{j-1}}=\left\lbrace\begin{array}{l}(\F{l_{j,\mu}}{\W\,\u\,.\,\mb{n_{j-1}}})\ \mr{on}\ \gamma_{j-1,1}\\
0 \hspace*{1.73cm} \mr{on}\ \gamma_{j-1,2}\end{array}\right.
\\z_{0,j,\mu|\gamma_{j-1}}= 0\end{array},\right.
\end{equation}
and, for $1\le k\le q$, $k\not=j$, $k\not=j-1$,
$$\left\lbrace\begin{array}{l}
(\di\F{\P z_{0,j,\mu}}{\P n})_{|\gamma_k}=0\\
z_{0,j,\mu|\gamma_k}= 0\end{array}.\hspace*{4.3cm}\right. $$
Next, applying Theorem \ref{4theu1}, let $z^*_{j,\mu} \in H^1(\O)$ be the unique solution of the problem
$$\left\lbrace\begin{array}{ll}
z^*_{j,\mu}+\W\,\u\,.\,\n z^*_{j,\mu}=l_{j,\mu}-z_{0,j,\mu}-\W\,\u\,.\,\n z_{0,j,\mu}\H &\mr{in}\H \O\\
z^*_{j,\mu}=0 &\mr{on}\H \G^-,\end{array}\right.$$
Then, $z_{j,\mu}=z^*_{j,\mu}+z_{0,j,\mu}$ verifies $z_{j,\mu}+\W\,\u\,.\,\n z_{j,\mu}=l_{j,\mu}$ and $z_{j,\mu|\G^-}=0$. Thus, in this first case, we have proven that
\begin{equation}\label{zjmu1}
z_{j,\mu}\in H^1(\O)\ \mr{is\ the\ solution\ of\ Problem}\ (P_{j,\mu}).\end{equation}
2) \tb{Second case}: $\mb{S^{-1}_j}\notin \gamma_{j-1}$, $\mb{S^{-1}_j}\notin E$, $\mb{S^{1}_j}\in \gamma_{j+1}$.\\
We can construct a lifting $z_{0,j,\mu}$ as in the first case. Since $l_{j,\mu}=0$ on $\gamma_k$ for $1\le k\le q$, $k\not=j$, there exists $z_{0,j,\mu}$ in $H^2(\O)$ verifying,  
$$\left\lbrace\begin{array}{l}
(\di\F{\P z_{0,j,\mu}}{\P n})_{|\gamma_j}=(\F{l_{j,\mu}}{\W\,\u\,.\,\mb{n_j}})_{|\gamma_j}\\
z_{0,j,\mu|\gamma_j}= 0\end{array}\right.\hspace*{2.15cm}$$
and, for $1\le k\le q$, $k\not=j$
$$\left\lbrace\begin{array}{l}
(\di\F{\P z_{0,j,\mu}}{\P n})_{|\gamma_k}=0\\
z_{0,j,\mu|\gamma_k}= 0\end{array}.\hspace*{4.3cm}\right. $$
Then, in the same way as in the first case, $z_{j,\mu}=z^*_{j,\mu}+z_{0,j,\mu}$ verifies $z_{j,\mu}+\W\,\u\,.\,\n z_{j,\mu}=l_{j,\mu}$ and  $z_{j,\mu|\G^-}=0$. Thus, in this second case, we have proven that
\begin{equation}\label{zjmu2}
z_{j,\mu}\in H^1(\O)\ \mr{is\ the\ solution\ of\ the\ problem}\ (P_{j,\mu}).\end{equation}
The cases where, for $i=-1,1$, $\mb{S^i_j}\in \gamma_{j+i}$ or $\mb{S^i_j}\notin \gamma_{j+i}$ with $\mb{S^i_j}\notin E$ can be studied in the same way as in the first two cases.\\
3) \tb{Third case}: $\mb{S^{-1}_j}\notin \gamma_{j-1}$, $\mb{S^{-1}_j}\in (E\cap S)$, $\mb{S^{1}_j}\in \gamma_{j+1}$.\\
\H Here, $[\mb{S^{-1}_j},\mb{S^{1}_j}]$ is a side of the polygon $\O$, $\mb{S^{-1}_j}$ is an end point of $\G^-$ such that $\u(\mb{S^{-1}_j})\,.\,\mb{n_j}=0$ and 
$\mb{S^{1}_j}$ is located inside $\G^-$. First, let us make the change of variables such that the point $\mb{S^{-1}_j}$ is the origin, the x-axis has the direction
of $\mb{n_j}$, oriented towards inside the domain $\O$, that is to say as the vector $-\mb{n_j}$, and with the segment $\gamma_j$ included in the positive $y$-axis, which is oriented by the tangent vector $\tho_-(\mb{S^{-1}_j})$(see the figure 3.6 below, where $\omega_j$ is the inner angle associated to the vertex $\mb{S^{-1}_j}$).\\
\unitlength=1cm
\begin{picture}(8,5.5)
\put(4,4){\vector(1,0){2.5}}
\put(1,4){\vector(0,-1){2.5}}
\put(0.5,1.5){$x$}
\put(2,1.5){Case $\omega_j\le \F{\pi}{2}$}
\put(6.1,3.7){$y$}
\qbezier(1.2,3.6)(1.35,3.7)(1.5,4)
\put(1.4,3.5){$\omega_j$}
\color{red}
\qbezier(4,4)(2.5,4)(1,4)
\put(3.8,4.3){$\mb{S^1_j}$}
\put(0.8,4.3){$\mb{S^{-1}_j}$}
\qbezier(4,4)(5,3)(6,2)
\put(2,4.3){$\gamma_j$}
\put(5.6,2.9){$\gamma_{j+1}$}
\color{black}
\put(3,4.3){\vector(0,1){0.7}}
\put(5.3,3.2){\vector(1,1){0.5}}
\put(3.1,4.3){$\vec{n}_j$}
\put(4.7,3.5){$\vec{n}_{j+1}$}
\put(2.2,3.75){\vector(1,0){0.7}}
\put(2.4,3.3){$\vec{\tho}_-(\mb{S^{-1}_j})$}
\put(7,1){figure\ 3.6}
\color{blue}
\qbezier(1,4)(1.5,3)(2,2)
\put(1.9,2.5){$\G^{0,+}$}
\color{green}
\put(1,4){\vector(-1,0){1}}
\put(-0.3,3.5){$\overrightarrow{u\scriptstyle ( S^{-1}_j)}$}
\end{picture}\hspace*{1cm}
\begin{picture}(8,5.5)
\put(4,4){\vector(1,0){2.5}}
\put(1,4){\vector(0,-1){2.5}}
\put(0.5,1.5){$x$}
\put(6.1,3.7){$y$}
\put(2,1.5){Case $\omega_j> \F{\pi}{2}$}
\qbezier(0.8,3.6)(1.2,3.65)(1.3,4)
\put(1.3,3.5){$\omega_j$}
\color{red}
\qbezier(4,4)(2.5,4)(1,4)
\put(3.8,4.3){$\mb{S^1_j}$}
\put(0.8,4.3){$\mb{S^{-1}_j}$}
\qbezier(4,4)(5,3)(6,2)
\put(2,4.3){$\gamma_j$}
\put(5.6,2.9){$\gamma_{j+1}$}
\color{black}
\put(3,4.3){\vector(0,1){0.7}}
\put(5.3,3.2){\vector(1,1){0.5}}
\put(3.1,4.3){$\vec{n}_j$}
\put(4.7,3.5){$\vec{n}_{j+1}$}
\put(2.2,3.75){\vector(1,0){0.7}}
\put(2.4,3.3){$\vec{\tho}_-(\mb{S^{-1}_j})$}
\color{blue}
\qbezier(1,4)(0.5,3)(0,2)
\put(-0.3,2.9){$\G^{0,+}$}
\color{green}
\put(1,4){\vector(-1,0){1}}
\put(-0.3,3.5){$\overrightarrow{u\scriptstyle ( S^{-1}_j)}$}
\end{picture}
\H With these new variables, since $\mb{S^{-1}_j}\in E$, we have 
$$\mb{S^{-1}_j}=(0,0),\ \u(\mb{S^{-1}_j})\,.\,\mb{n_{j}}=-u_1(0,0)=0$$ and the assumption (\ref{4ctheun}) yields $$\u(\mb{S^{-1}_j})\,.\,\tho_-(\mb{S^{-1}_j})=u_2(0,0)<0\ \mr{and}\ \di\F{\P\u}{\P\tho_-}(\mb{S^{-1}_j})\,.\,\mb{n_-(S^{-1}_j})=-\di\F{\P u_1}{\P y}(0,0)\not=0.$$ 
Considering that $\gamma_j\setminus\{\mb{S^{-1}_j}\}\subset\ov{\G^-}$, we have $u_1(0,y)>0$ for $y>0$ small enough. Thus, we derive the following properties of $\u$ in a neighborhood of $\mb{S^{-1}_j}=(0,0)$ :
\begin{equation}\label{4u1u2} 
u_1(0,0)=0,\ u_1(0,y)>0,\ u_2(0,0)<0\ \mr{and}\ \di\F{\P u_1}{\P y}(0,0)>0,\end{equation}
for $(0,y)\in \gamma_j\setminus \{\mb{S^{-1}_j}\}$, that is to say for $y>0$ small enough.\\
\H Next, we are going to split the problem $(P_{j,\mu})$ into two new problems. In this aim, we define a function $\lambda_{\mu}\in\mathcal{D}(\R^2)$ by
\begin{equation}\label{4lambdamu}
\forall\x\in\R^2,\ \lambda_{\mu}(\x)=\left\lbrace \begin{array}{l} 1\ \mr{if}
\ k(x,y)\le\mu\\[0.4cm]
0\ \mr{if}\ \ k(x,y)\ge 2\mu,\end{array}\right.\end{equation}
where $k(x,y)=|u_2(0,0)|\,|x|+\F{1}{2}|\F{\P u_1}{\P y}(0,0)|y^2$.\\
\H Then, we set 
\begin{equation}\label{4pbltl-}
\tilde{l}_{j,\mu}=(1-\lambda_{\mu})l_{j,\mu}\ \mr{and}\ \bar{l}_{j,\mu}=\lambda_{\mu}l_{j,\mu}\end{equation}
and we define the problem $(\tilde{P}_{j,\mu})$, which is associated to the right hand side $\tilde{l}_{j,\mu}$, and the problem 
$(\bar{P}_{j,\mu})$, which is associated to the right hand side $\bar{l}_{j,\mu}$. Since $l_{j,\mu}=\tilde{l}_{j,\mu}+\bar{l}_{j,\mu}$,  if we denote by, respectively, $z_{j,\mu}$, $\tilde{z}_{j,\mu}$ and $\bar{z}_{j,\mu}$ the unique solutions of, respectively, $(P_{j,\mu})$, $(\tilde{P}_{j,\mu})$ and $(\bar{P}_{j,\mu})$, we have
\begin{equation}\label{4dtpbp}
z_{j,\mu}= \tilde{z}_{j,\mu}+\bar{z}_{j,\mu}.\end{equation}
Thus, to prove that the problem $(P_{j,\mu})$ has its solution in $H^1(\O)$, we have only to prove that the problems $(\tilde{P}_{j,\mu})$ and $(\bar{P}_{j,\mu})$ have their solutions in $H^1(\O)$. Note that, extending the function $l\in H^1(\O)$ to $\R^2$, from now on, we will consider that the right hand sides $l$, $l_{j,\mu}$, $\tilde{l}_{j,\mu}$ and $\bar{l}_{j,\mu}$ belong to $H^1(\R^2)$  \\
\H First, we deal with the problem $(\tilde{P}_{j,\mu})$. Owing to the definition of the function $\lambda_{\mu}$, we can verify that $\tilde{l}_{j,\mu}$ vanishes on $\gamma_j$ on a neighborhood of the point $\mb{S^{-1}_j}$. So, we can construct a lifting $\tilde{z}_{0,j,\mu}$ in the same way as in the second case with $\tilde{l}_{j,\mu}$ in place of $l_{j,\mu}$, replacing $(\F{\tilde{l}_{j,\mu}}{\W\,\u\,.\,\mb{n_j}})_{|\gamma_j}$ with $0$ in a neighborhood of $\mb{S^{-1}_j}$ on $\gamma_j$ and $\tilde{z}_{j,\mu}=\z^*_{j,\mu}+\tilde{z}_{0,j,\mu}$ is the solution of the problem $(\tilde{P}_{j,\mu})$ in $H^1(\O)$.\\
\H Solving the problem $(\bar{P}_{j,\mu})$ is much more difficult because $\u(\mb{S^{-1}_j})\,.\,\mb{n_j}=0$ and $\bar{l}_{j,\mu}$ does not vanish in the neighborhood of $\mb{S^{-1}_j}$. From now on, we will use the following notation, for $r>0$ :
\begin{equation}\label{4dBmu} B_{j,r}=\{\x\in\R^2,\ d(\mb{S^{-1}_j},\x)=\sqrt{x^2+y^2}< r\},\ B_{j,r}^+=B_{j,r}\cap\{(x,y)\in\R^2,\ x\ge 0\}.\end{equation}
\H The proof will be built in several steps. In a first step, we define a local problem, which is the 
problem $(\bar{P}_{j,\mu})$ restricted to a neighborhhood $\O\cap B_{j,K}$ of $\mb{S^{-1}_j}$, and we express this local solution in integral form (see Lemma \ref{4zOK}). In a second step, we show that, if we chooze $\mu$ small enough, this local solution vanishes in $\ov{\O}\cap C(\mb{S^{-1}_j},r^*_{1,j},r^*_{2,j})$ where
$C(\mb{S^{-1}_j},r^*_{1,j},r^*_{2,j})$ is a ring centered in $\mb{S^{-1}_j}$ and included in $B_{j,K}$. In the third step, using its integral expression, we prove that the local solution belongs to 
$H^1(B_{j,r^*_{1,j}}\cap\O)$, which implies, owing to the second step, that its extension by zero is the solution $H^1$ of $(\bar{P}_{j,\mu})$.\\[0.2cm] 
\tb{First step}\\
\H In the following lemma, we give the expression of the local solution.
\begin{lem}\label{4zOK} Let $\mb{S^{-1}_j}$ belongs to $E\cap S$ and the real $K$ be defined by (\ref{4dK}). We set \linebreak$\O_{j,K}=\O\cap B_{j,K}$ and $\Gamma^-_{j,K}=\Gamma^-\cap B_{j,K}$. The solution of the problem
 $$\left\lbrace\begin{array}{ll}
z+\W\,\u\,.\,\n z=\bar{l}_{j,\mu}\H &\mr{in}\H \O_{j,K}\\
z=0 &\mr{on}\H \G^-_{j,K}\end{array}\right.$$
is expressed by 
\begin{equation}\label{4zsol}
\ z(x,y)=e^{-V(X(x,y),y)}(\int_{\a^{-1}(X(x,y))}^y \F{e^{V(X(x,y),t)}}{\scriptstyle\W U_2(X(x,y),t)}\bar{L}_{j,\mu}(X(x,y),t)\,dt),
\end{equation}
where $V$, $U_2$ and $\bar{L}_{j,\mu}$ are defined in (\ref{4dZLU2V}).
\end{lem}
\tb{Proof.} Owing to (\ref{4u1u2}), the continuity of $u_2$ yields that there exists a strictly positive real number $\mu_2\le \mu_0$, such that
\begin{equation}\label{4dmu2}
\forall\x=(x,y)\in B_{j,\mu_2}\cap\ov{\O},\ u_2(\x)<0.\end{equation}
In the same way, again the continuity of $u_2$ and the definition of $\di\F{\P u_1}{\P y}(0,0)$ with $u_1(0,0)=0$ imply that there exists a strictly positive real number $\mu_3\le \min(\mu_2,|\gamma_j|)$ such that
\begin{eqnarray}
\forall\x\in B_{j,\mu_3}\cap\ov{\O},\ \F{3}{2} u_2(0,0)\le u_2(\x)\le\F{1}{2} u_2(0,0)\ <0\ \ \mr{and}\nonumber\\ 
\forall y\in [0,\mu_3],\ \di\F{1}{2}\di\F{\P u_1}{\P y}(0,0)y\le u_1(0,y)\le \di\F{3}{2}\di\F{\P u_1}{\P y}(0,0)y \label{4dmu3}.\end{eqnarray}
\H For $0\le r_1<r_2$, let us 
define the sets
\begin{equation}\label{4dEr1r2}
E_{r_1,r_2}=\{(x,y)\in\R^2,\ r_1\le k(x,y)\le r_2\}\ \mr{and}\ C(\mb{S^{-1}_j},r_1,r_2)=\{(x,y)\in\R^2,\ r_1\le\sqrt{x^2+y^2}\le r_2\},\end{equation}
where $k(x,y)$ is defined in (\ref{4lambdamu}). 
Considering that, for $r>0$,
$$k(x,y)=r\Longleftrightarrow x^2+y^2=\F{(\F{\P u_1}{\P y}(0,0))^2}{4\,(u_2(0,0))^2}y^4+(\F{r}{|u_2(0,0)|})^2+(1-\F{r\,|\F{\P u_1}{\P y}(0,0)|}{(u_2(0,0))^2})y^2$$
and that $y^2\le \di\F{2r}{|\F{\P u_1}{\P y}(0,0)|}$, we can verify, for $0<r\le\di\F{(u_2(0,0))^2}{|\F{\P u_1}{\P y}(0,0)|}$, the following inclusions : 
\begin{equation}\label{4ikB}
B_{j,\F{r}{|u_2(0,0)|}}\subset E_{0,r}\subset B_{j,2\sqrt{\F{r}{|\F{\P u_1}{\P y}(0,0)|}}}\end{equation}
and if $(r_1,r_2)$ verifies 
\begin{equation}\label{4Cr1r2}
\left\lbrace\begin{array}{l}
0\le r_1<r_2\le\di\F{(u_2(0,0))^2}{|\F{\P u_1}{\P y}(0,0)|}\\
r_1<\di\F{|\F{\P u_1}{\P y}(0,0)|r_2^2}{4(u_2(0,0))^2}\end{array}\right.\ \mr{then}\ 
C(\mb{S^{-1}_j},2\sqrt{\F{r_1}{|\F{\P u_1}{\P y}(0,0)|}},\F{r_2}{|u_2(0,0)|})\subset E_{r_1,r_2}.
\end{equation}
\H Let us consider the transport equation $z+\W\,\u\,.\,\n z=\bar{l}_{j,\mu}$ of the problem $(\bar{P}_{j,\mu})$ and the following change of variables :
 we set for all $(x,y)\in \ov{\O}$ 
\begin{equation}\label{4chtv}
\left\lbrace\begin{array}{l}
X(x,y)=-\di\int_0^x u_2(t,y)\,dt+\di\int_0^y u_1(0,t)\,dt\\
Y(x,y)=y\end{array}\right..\end{equation}
\H Note that it is more convenient, especially in the case where $\o_j> \di\F{\pi}{2}$ ($\o_j$ is the inner angle of the polygon $\O$ associated to the vertex $\mb{S^{-1}_j}$), to define $X(x,y)$ 
for $y\le 0$ and $x\ge 0$ small enough. So, we will replace $u_2$ with an extension of $u_2$, defined for example by symmetries, in (\ref{4chtv}), this extension of $u_2$, still denoted $u_2$, verifying (\ref{4dmu2}), respectively (\ref{4dmu3}), in 
$B_{j,\mu_2}^+$, respectively $B_{j,\mu_3}^+$. Then, we define an extension of $u_1$, for $y\le 0$ and $x\ge 0$ small enough, by 
$$u_1(x,y)=-\int_0^x \F{\P u_2}{\P y}(t,y)\,dt,$$ such that $\d\u=0$ and $u_1(0,y)=0$, for $y\le 0$ small enough.
We can verify, in view of $\d\u=0$, that
\begin{equation}\label{4X'Y'}
X'_x=-u_2\ \mr{and}\ X'_y=u_1\ \mr{in}\ \ov{\O}.\end{equation} 
Let us show that the mapping 
\begin{equation}\label{4dphi}
\begin{array}{lll}\varphi:&B_{j,\mu_2}^+&\longrightarrow \varphi(B_{j,\mu_2}^+)\\
&(x,y)&\longmapsto (X,Y)\end{array}\ \mr{is\ one-to-one},\end{equation}
 where $\mu_2$ is defined in (\ref{4dmu2}). Let us assume that 
$$X(x,y)=X(x',y')\ \mr{and}\ Y(x,y)=Y(x',y')\ \mr{with}\ (x,y)\in B_{j,\mu_2}^+\ \mr{and}\ (x',y')\in B_{j,\mu_2}^+.$$ Then, the second equation gives directly $y=y'$ and we obtain $X(x,y)=X(x',y)$. Since $X'_x(x,y)=-u_2(x,y)>0$ for $(x,y)\in B_{j,\mu_2}^+$, we derive $x=x'$.\\
\H Since $\vfy$ is of class $C^1$ in $B_{j,\mu_2}^+$ and since the jacobian of the mapping $\vfy$ is $-u_2$, which is strictly positive in $B_{j\mu_2}^+$, we can define an inverse function 
$\vfy^{-1}$ of class $C^1$ in $\varphi(B_{j,\mu_2}^+)$. Then, in view of the definition of $\mu_3$ in (\ref{4dmu3}), we define the functions $Z$, $U_2$, $\bar{L}_{j,\mu}$ and $V$ on $\vfy(B_{j,\mu_3}^+)$ by
\begin{equation}\label{4dZLU2V}
 Z=z\circ \vfy^{-1}, \U_2=u_2\circ\vfy^{-1},\ \bar{L}_{j,\mu}=\bar{l}_{j,\mu}\circ\vfy^{-1}\ \mr{and}\ V:(X,Y)\mapsto \int_0^Y\F{1}{\W U_2(X,t)}\,dt.\end{equation}
\H Let us show that for $\x=(x,y)$ in a neighborhood of $(0,0)$ and $0\le |t|\le |Y|$ with $tY\ge 0$, then $(X(x,y),t)$ belongs to  
$\varphi(B_{j,\mu_2}^+)$. First, for $(x,y)\in B_{j,\F{\mu_3}{2}}^+$ and $0\le t\le Y=y$ or $y=Y\le t\le 0$ (case where $\omega_j>\F{\pi}{2}$), owing to (\ref{4dmu3}), in view of $u_1(0,y)=0$ for $y\le 0$, we have
$$X(0,t)\le X(x,y)\le \F{3}{2}k(x,y)\ \mr{and}\ X(\F{\mu_3}{2},t)\ge \F{\mu_3}{4}|u_2(0,0)|.$$  
Then, 
\begin{equation}\label{4dI}
\forall (x,y)\in \R^2\ \mr{such\ that}\H \left\lbrace\begin{array}{l}k(x,y\le \di\F{\mu_3|u_2(0,0)|}{6}\\[0.2cm]
(x,y)\in B_{j,\F{\mu_3}{2}}^+,\end{array}\right.\end{equation}
we have $X(0,t)\le X(x,y)\le X(\F{\mu_3}{2},t)$ and there exists a real number $x_t\in [0,\F{\mu_3}{2}]$ such that $X(x_t,t)=X(x,y))$ and, therefore, $(X(x,y),t)$ belongs to $\varphi(B_{j,\mu_3}^+)$. Finally, we set 
\begin{equation}\label{4dmu4tB}
\mu_4=\F{|u_2(0,0)|}{|\F{\P u_1}{\P y}(0,0)|}\ \mr{and}\ \widetilde{B}_j=B_{j,\min(\F{\mu_3}{6},\mu_4)}^+\cap\O.\end{equation} 
Since $\widetilde{B}_j=B_{j,\F{r}{|u_2(0,0)|}}^+\cap\O$ with $r=\min(\di\F{\mu_3|u_2(0,0)|}{6},\di\F{(u_2(0,0))^2}
{|\F{\P u_1}{\P y}(0,0)|})$, 
in view of (\ref{4ikB}), all $\x=(x,y)\in \widetilde{B}_j$ verifies (\ref{4dI}) and, consequently, $(X(x,y),t)$ belongs to $\varphi(B_{j,\mu_3}^+)$. Then, with the new functions defined in (\ref{4dZLU2V}), in view of $\mr{div}\,\u=0$, we have the following equivalence :
 $$z+\W\,\u\,.\,\n z=\bar{l}_{j,\mu}\ a.e.\ \mr{in}\ \widetilde{B}_j\iff Z+\W\,U_2\F{\P Z}{\P Y}=\bar{L}_{j,\mu}\ a.e.
 \ \mr{in}\ \varphi(\widetilde{B}_j).$$ 
 Solving this last equation yields
 \begin{eqnarray*}\forall(X,Y)\in\varphi(\widetilde{B}_j),\ Z(X,Y)=e^{-V(X,Y)}(\int_0^Y \F{e^{V(X,t)}}{\scriptstyle\W U_2(X,t)}\bar{L}_{j,\mu}(X,t)\,dt+C(X))\Longleftrightarrow\\
 \forall (x,y)\in\tilde{B}_j,\ z(x,y)=e^{-V(X(x,y),y)}(\int_0^y\F{e^{V(X(x,y),t)}}{\mathcal{W} U_2(X(x,y),t)}\bar{L}_{j,\mu}(X(x,y),t)\,dt+C(X(x,y))),\end{eqnarray*}
 where $C$ is a function of $L^2$. We have to compute the function $C$ so that the solution $Z$ verifies the boundary condition on $\G^-$.\\ 
\H Let us define the real number $y_M>0$  by
$$ y_M=\sup\{y,\ \m(0,y)\in \widetilde{B}_j\cap\gamma_j\}$$
and the function $\a$ on the set $[0,y_M]$ by
\begin{equation}\label{4defaa-1}
\forall y\in [0,y_M],\ \a(y)=X(0,y).\end{equation}
Note that \begin{equation}\label{4yM}
y_M=\min(\F{\mu_3}{6},\mu_4,|\gamma_j|).\end{equation}
Considering that, $\forall y\in]0,y_M]$,$$ \a'(y)=u_1(0,y)>0,$$ 
the mapping $\a$ from $[0,y_M]$ to $[0,\a(y_M)]$ is one-to-one and we can define the inverse function $\a^{-1}$ from $[0,\a(y_M)]$ to $[0,y_M]$. 
Moreover, $\a^{-1}$ is strictly positive on $]0,\a(y_M)]$. Then, the continuity of the functions $X$ and 
$Y$ yields that there exist a real number $\mu_5>0$ such that
\begin{equation}\label{4cXY}
\forall (x,y)\in B_{j,\mu_5}^+,\ X(x,y)\in [0,\a(y_M)].\end{equation}
Finally, we set
\begin{equation}\label{4dK}
K=\min(\F{\mu_3}{6},\mu_4,\mu_5,|\gamma_j|),\end{equation}
where the constants $\mu_3$, $\mu_4$ and $\mu_5$ are defined, respectively, by (\ref{4dmu3}), (\ref{4dmu4tB}) and (\ref{4cXY}).\\
Then, the boundary condition $z_{|\gamma_j}=0$ allows us to compute the function $C$. Indeed, setting $s=X(0,y)=\a(y)\Leftrightarrow y=\a^{-1}(s)$, we have
\begin{eqnarray*}z_{|\G^-_{j,K}}=0\Leftrightarrow \forall y,\ 0\le y\le K,\ z(0,y)=0
\Leftrightarrow \forall s,\ 0\le s\le \a(K),\ Z(s,\a^{-1}(s))=0\\ \Leftrightarrow\forall s,\ 0\le s\le \a(K),\ C(s)=-\int_0^{\a^{-1}(s)}\F{e^{V(s,t)}}{\mathcal{W}U_2(s,t)}\bar{L}_{j,\mu}(s,t)\,dt\hspace*{1.6cm}\end{eqnarray*}
and we obtain that the solution $z$ is expressed in $\O_{j,K}$ as (\ref{4zsol}).
\hfill$\diamondsuit$\\[0.2cm]
\tb{Second step}\\
\H Let us show that, for $\mu$ small enough and $\x$ far enough from $\mb{S^{-1}_j}$, then  $z(x,y)=0$. More precisely, let us prove the following lemma.
\begin{lem}\label{4lz=0}
Let $r_{1,j}$ and $r_{2,j}$ be defined by (\ref{4dr1r2}), let $\mu>0$ such that $\mu\le \F{r_{1,j}}{6}$ and let the local solution $z$ of the problem $(\bar{P}_{j,\mu})$ be expressed by (\ref{4zsol}). Then, 
\begin{equation}\label{4z=0}
\forall(x,y)\in C(\mb{S^{-1}_j},2\sqrt{\F{r_{1,j}}{|\F{\P u_1}{\P y}(0,0)|}},\F{r_{2,j}}{|u_2(0,0)|})\cap\ov{\O},\ z(x,y)=0.\end{equation}
\end{lem}
\tb{Proof.} Let us note that, if $(x,y)\in B^+_{j,K}$, then $y\le \a^{-1}(X(x,y))$. Indeed, if $y<0$, then $y<0\le \a^{-1}(X(x,y))$ and if $y\ge 0$, then 
$\a(y_M)\ge X(x,y)\ge X(0,y)=\a(y)\ge 0$, which implies $y\le \a^{-1}(X(x,y))$, since $\a^{-1}$ is strictly increasing on $[0,\a(y_M)]$. Thus, we distinguish two cases :\\
a) \un{First case} : $(x,y)\in B_{j,\F{K}{6}}^+$ and $0\le y\le t\le\a^{-1}(X(x,y))$.\\
Note that, since $\F{K}{6}\le \mu_5$, then $\a^{-1}(X(x,y))\le y_M$, which implies $t\le\F{\mu_3}{6}$.
On the one hand, we have $$\a(t)=X(0,t)\le X(x,y).$$ 
On the other hand, since $|x|\le \F{K}{6}$, we derive
$$X(x,y)\le \F{3}{2}|u_2(0,0)|\,|x|+\int_0^y u_1(0,\theta)\,d\theta\le \F{1}{4}|u_2(0,0)|\,K+\int_0^t u_1(0,\theta)\,d\theta\le X(\F{K}{2},t).$$
Therefore, if $(x,y)\in B_{j,\F{K}{6}}^+$ with $y\ge 0$, there exists $x_t\in [0,\F{K}{2}]$ such that
\begin{equation}\label{4dtxbis}
 X(x,y)=X(x_t,t)\ \mr{with}\ (x_t,t)\in B_{j,\mu_3}^+.\end{equation}
Then, the inequalities (\ref{4dmu3}) yield 
\begin{equation}\label{4ixtx}
\F{1}{2}k(x,y)\le X(x,y)=X(x_t,t)\le \F{3}{2}k(x_t,t)\Longrightarrow\F{1}{3}k(x,y)\le k(x_t,t).
\end{equation} 
We set
\begin{equation}\label{4dr1r2}
r_{1,j}=\min(\F{|u_2(0,0)|K}{12},\F{|\F{\P u_1}{\P y}(0,0)| K^2}{288})\ \mr{and}\ r_{2,j}=\F{K|u_2(0,0)|}{6}.\end{equation}
Choosing the real number $\mu>0$ such that
\begin{equation}\label{4cmu}
\mu\le \F{r_{1,j}}{6}\Longleftrightarrow 6\mu\le r_{1,j},\end{equation} 
we can verify that $0<r_{1,j}<r_{2,j}\le \di\F{(u_2(0,0))^2}{|\F{\P u_1}{\P y}(0,0)|}$ and $r_{1,j}<\di\F{|\F{\P u_1}{\P y}(0,0)|r_{2,j}^2}{4(u_2(0,0))^2}$ and, owing to 
(\ref{4Cr1r2}), we obtain, 
$$\forall(x,y)\in C(\mb{S^{-1}_j},2\sqrt{\F{r_{1,j}}{|\F{\P u_1}{\P y}(0,0)|}},\F{r_{2,j}}{|u_2(0,0)|})\cap\ov{\O}\ \mr{with}\ y\ge 0,\ r_{1,j}\le k(x,y)\le r_{2,j}.$$
Hence, in view of (\ref{4ixtx}), we derive
\begin{eqnarray*}
\forall(x,y)\in C(\mb{S^{-1}_j},2\sqrt{\F{r_{1,j}}{|\F{\P u_1}{\P y}(0,0)|}},\F{r_{2,j}}{|u_2(0,0)|})\cap\ov{\O}\ \mr{with}\ y\ge 0,\\
 k(x_t,t)\ge \F{1}{3}k(x,y)\ge\F{1}{3}r_{1,j}\ge 2\mu.\end{eqnarray*}
Finally, since $\bar{L}_{j,\mu}(X(x,y),t)=\bar{L}_{j,\mu}(X(x_t,t),t)=\bar{l}_{j,\mu}(x_t,t)$, considering (\ref{4lambdamu}), (\ref{4pbltl-}) and (\ref{4zsol}), 
we obtain 
$$
\forall(x,y)\in C(\mb{S^{-1}_j},2\sqrt{\F{r_{1,j}}{|\F{\P u_1}{\P y}(0,0)|}},\F{r_{2,j}}{|u_2(0,0)|})\cap\ov{\O}\ \mr{with}\ y\ge 0,\ z(x,y)=0,
$$
where $r_{1,j}$, $r_{2,j}$ and $\mu$ are given by (\ref{4dr1r2}) and (\ref{4cmu}).\\
b) \un{Second case} ($\omega_j>\di\F{\pi}{2}$) : $(x,y)\in B_{j,\F{K}{6}}^+$, $y<0$ and $y\le t\le\a^{-1}(X(x,y))$.\\
\H Choosing first $t\in [0, \a^{-1}(X(x,y))]$ and second $t\in[y,0]$, considering that $u_1(0,y)=0$ when $y<0$,  we process in the same way as previously
and we obtain, as in the case where $y\ge 0$, that there exists $x_t\in [0,\F{K}{2}]$ such that
\begin{equation}\label{4dtxbis-}
 X(x,y)=X(x_t,t)\ \mr{with}\ (x_t,t)\in B_{j,\mu_3}^+\end{equation}
and $\forall(x,y)\in C(\mb{S^{-1}_j},2\sqrt{\F{r_{1,j}}{|\F{\P u_1}{\P y}(0,0)|}},\F{r_{2,j}}{|u_2(0,0)|})\cap\ov{\O}$ with $y< 0,$
$$\int_{\a^{-1}(X(x,y))}^0 \F{e^{V(X(x,y),t)}}{\scriptstyle\W U_2(X(x,y),t)}\bar{L}_{j,\mu}(X(x,y),t)\,dt=0,\ \int_0^y \F{e^{V(X(x,y),t)}}{\scriptstyle\W U_2(X(x,y),t)}\bar{L}_{j,\mu}(X(x,y),t)\,dt=0,$$
which implies $z(x,y)=0$. Finally, gathering the cases $y\ge 0$ and $y<0$, we derive (\ref{4z=0}).\hfill$\diamondsuit$\\[0.2cm]
\tb{Third step}\\
\H Next, we will prove the following lemma that gives the regularity $H^1$ of the local solution of the problem $(\bar{P}_{j,\mu})$.
\begin{lem} Let $z$ be defined by (\ref{4zsol})and let $r^*_j$ be defined by 
\begin{equation}\label{4dr*j}
r^*_j=2\sqrt{\F{r_{1,j}}{|\F{\P u_1}{\P y}(0,0)|}}=\min( \sqrt{\F{|u_2(0,0)K}{3|\F{\P u_1}{\P y}(0,0)|}},\F{K}{6\sqrt{2}}),\end{equation}
where $r_{1,j}$ and $K$ are defined in (\ref{4dr1r2}) and (\ref{4dK}).
Then $z$ belongs to $H^1(B_{j,r^*_j}\cap\O)$.\end{lem}
\tb{Proof.} Let us prove first that $z$ belongs to $L^2(B_{j,r^*_j}\cap\O)$. Using the change of variables defined in (\ref{4chtv}) yields $z(x,y)=Z(X,Y)$ with
\begin{equation}\label{4Z(X,Y)}Z(X,Y)=e^{-V(X,Y)}\int_{\a^{-1}(X)}^Y\F{e^{V(X,t)}}{\mathcal{W} U_2(X,t)}\bar{L}_{j,\mu}(X,t)\,dt\end{equation} and, in view of the jacobian $\F{D(x,y)}{D(X,Y)}=-\F{1}{U_2}$ and the inequality of Cauchy-Schwarz, we obtain
\begin{eqnarray}
\int\int_{B_{j,r^*_j}\cap\O}z^2\,dxdy=\int\int_{\varphi(B_{j,r^*_j}\cap\O)}Z^2\F{1}{|U_2|}\,dXdY\hspace*{5cm}\label{4Zi}\\
\le \int\int_{\varphi(B_{j,r^*_j}\cap\O)}\F{e^{2V}}{|U_2|}(\int_Y^{\a^{-1}(X)}\F{e^{-2 V(X,t)}}{\mathcal{W}^2(U_2(X,t))^2}\,dt)
(\int_Y^{\a^{-1}(X)}(\bar{L}_{j,\mu}(X,t))^2\,dt)dXdY.\nonumber\end{eqnarray}
We have to estimate the terms of the previous integral. Since $r^*_j\le\F{K}{6}$, in view of (\ref{4dmu3}), we derive
$$ \forall (x,y)\in B_{j,r^*_j}\cap\O,\ \F{1}{|U_2(X,Y)|}=\F{1}{|u_2(x,y))|}\le \F{2}{|u_2(0,0)|}.$$
Owing to (\ref{4dmu4tB}) and $(B_{j,r^*_j}\cap\O)\subset \widetilde{B_j}$, for $(X,Y)\in\varphi(B_{j,r^*_j}\cap\O)$ and  $|t|\le |Y|$, we have $U_2(X,t)=u_2(x_t,t)$ with $x_t\in [0,\F{\mu_3}{2}]$ and, considering that $|Y|\le \F{K}{6}$, we obtain
$$|V(X,Y)|=|\int_0^Y\F{1}{\mathcal{W}u_2(x_t,t)}\,dt\le \F{K}{3|\mathcal{W}||u_2(0,0)|}.$$
In the same way, for $(X,Y)\in\varphi(B_{j,r^*_j}\cap\O)$ and $Y\le t\le\a^{-1}(X)$, in view of (\ref{4dtxbis}) and (\ref{4dtxbis-}), we have 
$U_2(X,t)=U_2(X(x_t,t),t)=u_2(x_t,t)$ with $(x_t,t)\in B_{j,\mu_3}^+$, which implies
$$\F{1}{|U_2(X,t)|}\le \F{2}{|u_2(0,0)|}$$
and, for $V(X,t)=\di\int_0^t \F{1}{\mathcal{W}U_2(X,\theta)}\,d\theta$, we prove that $X=X(x_{\theta},\theta)$ with $(x_{\theta},\theta)\in B_{j,\mu_3}^+$, which gives, since $-\F{K}{6}\le Y\le t\le \a^{-1}(X)\le y_M\le \F{\mu_3}{6}$, the following estimate
$$|V(X,t)|\le \F{\mu_3}{3|\mathcal{W}||u_2(0,0)|}.$$
Substituting these bounds in (\ref{4Zi}) yields that there exists a strictly positive constant $C$ such that 
\begin{equation}\label{4miz}
\int\int_{B_{j,r^*_j}\cap\O} z^2\,dxdy\le C\int\int_{\varphi(B_{j,r^*_j}\cap\O)}(\int_Y^{\a^{-1}(X)}(\bar{L}_{j,\mu}(X,t))^2\,dt)dXdY.\end{equation}
Since $\forall (X,Y)\in\varphi(B_{j,r^*_j}\cap\O)$, we have $-\F{K}{6}\le Y\le\F{K}{6}$ and $X(0,Y)\le X\le X(\F{K}{6},Y)$, we obtain
\begin{eqnarray*}
\int\int_{B_{j,r^*_j}\cap\O }z^2\,dxdy\le C\int_{-\F{K}{6}}^{\F{K}{6}}dY\left(\int_{X(0,Y)}^{X(\F{K}{6},Y)}\,dX
(\int_Y^{\a^{-1}(X)}(\bar{L}_{j,\mu}(X,t)))^2\,dt)\right),\\
\le C\int_{-\F{K}{6}}^{\F{K}{6}}dY\left (\int\int_{D_Y}(\bar{L}_{j,\mu}(X,t)))^2\,dX\,dt\right),\end{eqnarray*}
where $D_Y=\{(X,t)\in\R^2,\ X(0,Y)\le X\le X(\F{K}{6},Y),\ Y\le t\le \a^{-1}(X)\}$. Next, we compute the integral on $D_Y$ by making the substitution 
$\left\lbrace\begin{array}{l} X=X(\tilde{x},t)\\
t=t\end{array}\right.$, the jacobian of which is $-u_2(\tilde{x},t)$. Indeed, the mapping 
$\begin{array}{lll}\psi:&D_Y&\longrightarrow \psi(D_Y)\\
&(X,t)&\longmapsto (\tilde{x},t)\end{array}\ \mr{is\ one-to-one}$ and of class $C^1$ on $D_Y$, as we proved previously by (\ref{4dtxbis}) and (\ref{4dtxbis-}),
with $\tilde{x}\in[0,\F{K}{2}]$ and $t\in[-\F{K}{6},\F{\mu_3}{6}]$. Thus, the jacobian is strictly positive and bounded by $\F{3}{2}|u_2(0,0)|$, and we obtain
\begin{equation}\label{4zL2}
\int\int_{B_{j,r^*_j}\cap\O} z^2\,dxdy\le \F{1}{2}CK|u_2(0,0)|\int\int_{[0,\F{K}{2}]\t [-\F{K}{6},\F{\mu_3}{6}]}\bar{l}_{j,\mu}(\tilde{x},t)\,d\tilde{x}\,dt<+\iy,\end{equation}
which proves that $z$ belongs to $L^2(B_{j,r^*_j}\cap\O)$, since $\bar{l}_{j,\mu}$ belongs to $H^1(\R^2)$.\\
\H It remains to prove that $\n z$ belongs to $L^2(B_{j,r^*_j}\cap\O)$. Again, we use the change of variables defined in (\ref{4chtv}). Computing the partial derivatives yields
$$\F{\P z}{\P x}=-U_2\F{\P Z}{\P X}\ \mr{and}\ \F{\P z}{\P y}=U_1\F{\P Z}{\P X}+\F{\P Z}{\P Y}.$$
Then, the inequality $2 |ab|\le a^2+b^2$ implies
$|\n z|^2\le (1+U_1^2+U_2^2)|\n Z|^2$, where $|\ .\ |$ represents the euclidian norm in $\R^2$. Hence, we derive
$$\int\int_{B_{j,r^*_j}\cap\O}|\n z|^2\,dxdy\le \F{2\max\limits_{x\in\ov{\O}}(1+u_1^2+u_2^2)}{|u_2(0,0)|}\int\int_{\varphi(B_{j,r^*_j}\cap\O)}|\n Z|^2\,dXdY.$$ 
Since $Z+\mathcal{W}U_2\di\F{\P Z}{\P Y}=\bar{L}_{j,\mu}$, owing to (\ref{4dmu3}), we obtain that $\di\F{\P Z}{\P Y}$ belongs to $L^2(\varphi(B_{j,r^*_j}\cap\O))$.\\
\H Next, we now come to the crucial point, which is to prove that the other partial derivative 
$\di\F{\P Z}{\P X}$ belongs to $L^2(\varphi(B_{j,r^*_j}\cap\O))$. From (\ref{4Z(X,Y)}), computing this derivative yields
\begin{eqnarray}\F{\P Z}{\P X}(X,Y)=-\F{\P V}{\P X}(X,Y)Z(X,Y)+e^{-V(X,Y)}\hspace*{6.3cm}\label{4dZdX}\\(\int_{\a^{-1}(X)}^Y \F{\P(\F{e^{V(X,t)}}{\mathcal{W}U_2(X,t)}\bar{L}_{j,\mu}(X,t))}{\P X}\,dt-
(\a^{-1})'(X)(\F{e^{V(X,\a^{-1}(X))}}{\mathcal{W}U_2(X,\a^{-1}(X))}\bar{L}_{j,\mu}(X,\a^{-1}(X)))\nonumber\end{eqnarray}
and
$$\F{\P(\F{e^{V(X,t)}}{\mathcal{W}U_2(X,t)}\bar{L}_{j,\mu}(X,t))}{\P X}=e^{V(X,t)}\F{\F{\P V}{\P X}(X,t)U_2(X,t)-\F{\P U_2}{\P X}(X,t)}{(U_2(X,t))^2}
 \bar{L}_{j,\mu}(X,t)+\F{e^{V(X,t)}}{U_2(X,t)}\F{\P\bar{L}_{j,\mu}}{\P X}(X,t).$$
 Then, for $(X,Y)\in\varphi(B_{j,r^*_j}\cap\O)$ and $|t|\le |Y|$ or $Y\le t\le\a^{-1}(X)$, we prove that $(X,t)$ belongs to $\varphi(B_{j,\mu_3}^+)$. Considering that $U_2$ is strictly negative and of class $C^1$ on $\varphi(B_{j,\mu_3}^+)$, we derive that $U_2$ and 
 $V$ are of class $C^1$ on $\varphi(B_{j,\mu_3}^+)$, which implies that the functions $\F{\P U_2}{\P X}$ and $\F{\P V}{\P X}$ are bounded on 
 $\varphi(B_{j,\mu_3}^+)$. Hence, there exist strictly positive constants $C_1$, $C_2$, $C_3$ and $C_4$ such that
 \begin{eqnarray}
 |\F{\P Z}{\P X}(X,Y)|\le C_1 |Z(X,Y)|+C_2\int_Y^{\a^{-1}(X)} |\bar{L}_{j,\mu}(X,t)|\,dt\nonumber\hspace*{3.8cm}\\
 +C_3\int_Y^{\a^{-1}(X)} |\F{\P\bar{L}_{j,\mu}}{\P X}(X,t)|\,dt+C_4
 |(\a^{-1})'(X)||\bar{L}_{j,\mu}(X,\a^{-1}(X))|.\label{4midZdX}\end{eqnarray}
 In view of the equalities $(X,\a^{-1}(X))= (\a(\a^{-1}(X)),\a^{-1}(X)=(X(0,\a^{-1}(X)),\a^{-1}(X))$, we derive $\bar{L}_{j,\mu}(X,\a^{-1}(X))=\bar{l}_{j,\mu}(0,\a^{-1}(X))$. Next, owing to $|\a^{-1}(X)-Y|\le \F{\mu_3+K}{6}$, using inequalities of Cauchy-Schwarz and 
 setting $C_5=C_1^2+C_2^2+C_3^2+C_4^2$ and $C_6=\F{\mu_3+K}{6}C_5$ yield 
 \begin{eqnarray} 
(\F{\P Z}{\P X}(X,Y))^2\le C_5((Z(X,Y))^2+((\a^{-1})'(X))^2(\bar{l}_{j,\mu}(0,\a^{-1}(X)))^2)\nonumber\\
+C_6\int_Y^{\a^{-1}(X)}((\bar{L}_{j,\mu}(X,t))^2+(\F{\P\bar{L}_{j,\mu}}{\P X}(X,t))^2)\,dt.\hspace*{0.6cm}\label{4midZdX2}\end{eqnarray}
There is only one term that is difficult to bound in $L^2(\varphi(B(\mb{S^{-1}_j},r^*_j)\cap\ov{\O}))$. Indeed, we have just prove that $Z$ 
belongs to $L^2(\varphi(B_{j,r^*_j}\cap\O))$ and for $$\int_Y^{\a^{-1}(X)}((\bar{L}_{j,\mu}(X,t))^2+(\F{\P\bar{L}_{j,\mu}}{\P X}(X,t))^2)\,dt,$$ we apply the previous method, which allowed us to bound the right hand side of the inequality (\ref{4miz}), using the same substitution $\left\lbrace\begin{array}{l} X=X(\tilde{x},t)\\
t=t\end{array}\right.$ as previously, since $\bar{l}_{j,\mu}$ belongs to 
$H^1(\R^2)$ and since we have $\di \F{\P \bar{L}_{j,\mu}}{\P X}(X,t)=-\di\F{1}{u_2(\tilde{x},t)}\di\F{\P \bar{l}_{j,\mu}}{\P \tilde{x}}(\tilde{x},t)$. It remains to bound the basic term 
$$\int\int_{\varphi(B_{j,r^*_j}\cap\O)}((\a^{-1})'(X))^2(\bar{l}_{j,\mu}(0,\a^{-1}(X)))^2\, dXdY.$$
Let us recall that, in view of (\ref{4defaa-1}),(\ref{4cXY}), since $(x,y)$ belongs to $B_{j,r^*_j}\cap\ov{\O}$, we have 
$$X\in [0,\a(y_M)]\Longleftrightarrow \a^{-1}(X)\le y_M=\le\min(\F{\mu_3}{6},\mu_4,|\gamma_j|).$$
 Since $\a'(\a^{-1}(X))=u_1(0,\a^{-1}(X))>0$ and considering (\ref{4dmu3}), on the one hand, we derive
\begin{equation}\label{4ma'-1}
|(\a^{-1})'(X)|=\F{1}{|\a'(\a^{-1}(X))|}=\F{1}{|u_1(0,\a^{-1}(X)|}\le \F{2}{|\F{\P u_1}{\P y}(0,0)||\a^{-1}(X)|}.\end{equation}
On the other hand, owing again to (\ref{4dmu3}), we have
$$|X|=|\a(\a^{-1}(X)|=|\int_0^{\a^{-1}(X)}u_1(0,\theta)\,d\theta|\le \F{3}{4}
|\F{\P u_1}{\P y}(0,0)|(\a^{-1}(X))^2.$$Substituting this inequality in (\ref{4ma'-1}) yields the following basic estimate of $|(\a^{-1})'(X)|$
\begin{equation}\label{4ma'-12}
|(\a^{-1})'(X)|\le \sqrt{\F{3}{|\F{\P u_1}{\P y}(0,0)|}}\F{1}{\sqrt{X}}.\end{equation}
Next, we distinguish two cases : if the angle $\omega_j\le \F{\pi}{2}$, then
$$\varphi(B_{j,r^*_j}\cap\O)\subset E_1=\{(X,Y)\in\R^2,\ Y\in [0,\F{K}{6}],\ X\in [X(0,Y),X(\F{K}{6},Y)],$$\pagebreak\\
and, if the angle $\omega_j> \F{\pi}{2}$, then
$$\varphi(B_{j,r^*_j}\cap\O)\subset E_1\cup E_2,$$
where 
$$ E_2=\{(X,Y)\in\R^2,\ Y\in [-\F{K}{6},0],\ X\in [X((\tan\omega_j)Y,0),X(\F{K}{6},0)].$$
Therefore, in the both case, we have to compute, for $Y\ge 0$
$$\int_{X(0,Y)}^{X(\F{K}{6},Y)}((\a^{-1})'(X))^2(\bar{l}_{j,\mu}(0,\a^{-1}(X)))^2\, dX.$$
Considering that $\bar{l}_{j,\mu}$ belongs to $H^1(\R^2)$, we derive that the function $(0,y)\mapsto \bar{l}_{j,\mu}(0,y)$ belongs to 
$H^{\F{1}{2}}(\gamma_j)\supset L^6(\gamma_j)$. Hence, in view of (\ref{4ma'-12}), using the Holder's inequality 
 yields
\begin{eqnarray*}
\int_{X(0,Y)}^{X(\F{K}{6},Y)}((\a^{-1})'(X))^2(\bar{l}_{j,\mu}(0,\a^{-1}(X)))^2\, dX\hspace*{4cm}\\ 
\le\F{3}{|\F{\P u_1}{\P y}(0,0)|}(\int_{X(0,Y)}^{X(\F{K}{6},Y)}\F{1}{X^{\F{3}{2}}}dX)^{\F{2}{3}}(\int_{X(0,Y)}^{X(\F{K}{6},Y)} (\bar{l}_{j,\mu}(0,\a^{-1}(X)))^6\, dX)^{\F{1}{3}}.\end{eqnarray*}
Setting $v=\a^{-1}(X)$, with (\ref{4dmu3}), we obtain
$$\int_{X(0,Y)}^{X(\F{K}{6},Y)} (\bar{l}_{j,\mu}(0,\a^{-1}(X)))^6\,dX\le\int_0^{\a(y_M)}(\bar{l}_{j,\mu}(0,v))^6|u_1(0,v)|\,dv
\le \F{3}{2}\a(y_M)|\F{\P u_1}{\P y}(0,0)|\,\p\bar{l}_{j,\mu}\p^6_{L^6(\gamma_j)},$$which implies
$$\int_{X(0,Y)}^{X(\F{K}{6},Y)}((\a^{-1})'(X))^2(\bar{l}_{j,\mu}(0,\a^{-1}(X)))^2\, dX
\le (\F{162\a(y_M)}
{(\F{\P u_1}{\P y}(0,0))^2})^{\F{1}{3}}\F{\p\bar{l}_{j,\mu}\p^2_{L^6(\gamma_j)}}{(X(0,Y))^{\F{1}{3}}}.$$
Hence,  since $\di\F{1}{|X(0,Y)|}\le \di\F{4}{|\F{\P u_1}{\P y}(0,0)|Y^2}$ (see (\ref{4dmu3})), we derive\begin{equation}\label{4mInt1}
\int_{X(0,Y)}^{X(\F{K}{6},Y)}((\a^{-1})'(X))^2(\bar{l}_{j,\mu}(0,\a^{-1}(X)))^2\, dX\le \F{(648\a(y_M))^{\F{1}{3}}
\p\bar{l}_{j,\mu}\p^2_{L^6(\gamma_j)}}{|\F{\P u_1}{\P y}(0,0)|Y^{\F{2}{3}}}.\end{equation}
In the case where $\omega_j> \F{\pi}{2}$, we have to bound, for $Y<0$, 
$$\int_{X((\tan\omega_j)Y,0)}^{X(\F{K}{6},0}((\a^{-1})'(X))^2(\bar{l}_{j,\mu}(0,\a^{-1}(X)))^2\, dX.$$
In the same way as previously, we obtain
$$\int_{X((\tan\omega_j)Y,0)}^{X(\F{K}{6},0)}((\a^{-1})'(X))^2(\bar{l}_{j,\mu}(0,\a^{-1}(X)))^2\, dX
\le (\F{162\a(y_M)}
{(\F{\P u_1}{\P y}(0,0))^2})^{\F{1}{3}}\F{\p\bar{l}_{j,\mu}\p^2_{L^6(\gamma_j)}}{(X((\tan\omega_j)Y,0)))^{\F{1}{3}}}.$$
Since $|X((\tan\omega_j)Y,0)|\ge \di\F{2}{|u_2(0,0|\,|\tan\omega_j|\,|Y|}$, we derive
\begin{equation}\label{4mInt2}
\int_{X((\tan\omega_j)Y,0)}^{X(\F{K}{6},0)}((\a^{-1})'(X))^2(\bar{l}_{j,\mu}(0,\a^{-1}(X)))^2\, dX\le\F{(324\a(y_M))^{\F{1}{3}}\p\bar{l}_{j,\mu}\p^2_{L^6(\gamma_j)}}{((\F{\P u_1}{\P y}(0,0))^2|u_2(0,0)|\,|\tan\omega_j|)^{\F{1}{3}}}\F{1}{|Y|^{\F{1}{3}}}.\hspace*{0.8cm}\end{equation}
In the case where $\omega_j\le \F{\pi}{2}$, integrating with respect to $Y$ on the interval $[0,\F{K}{6}]$ the both side of (\ref{4mInt1}) yields
 $$\int\int_{\varphi(B_{j,r^*_j}\cap\ov{\O})}((\a^{-1})'(X))^2(\bar{l}_{j,\mu}(0,\a^{-1}(X)))^2\, dXdY\le 9\F{(4\a(y_M)K)^{\F{1}{3}}}{|\F{\P u_1}{\P y}(0,0))|}\p\bar{l}_{j,\mu}\p^2_{L^6(\gamma_j)}<+\iy.$$
 In the case where $\omega_j> \F{\pi}{2}$, in addition to the previous integral, we must integrate with respect to $Y$ on the interval $[-\F{K}{6},0]]$ the both side of (\ref{4mInt2}), which gives\begin{eqnarray*}
 \int\int_{\varphi(B_{j,r^*_j}\cap\ov{\O})}((\a^{-1})'(X))^2(\bar{l}_{j,\mu}(0,\a^{-1}(X)))^2\, dXdY\hspace*{3.2cm}\\
 \le (9\F{(4\a(y_M)K)^{\F{1}{3}}}{|\F{\P u_1}{\P y}(0,0))|}+ \F{3}{2}\F{(9\a(y_M)K^2)^{\F{1}{3}}}{((\F{\P u_1}{\P y}(0,0))^2|u_2(0,0)|\,|\tan\omega_j|)^{\F{1}{3}}})\p\bar{l}_{j,\mu}\p^2_{L^6(\gamma_j)}<+\iy.\end{eqnarray*}
Finally, in view of (\ref{4midZdX2}), we have obtained that $\di\F{\P Z}{\P X}$ belongs to $L^2(\varphi(B_{j,r^*_j}\cap\O))$, which implies, as we saw previously, that 
$\n z$ belongs to $L^2(B_{j,r^*_j}\cap\O)$ and, therefore,
 with (\ref{4zL2}), we derive that $z$ belongs to $H^1( B_{j,r^*_j}\cap\O)$, which ends the proof of the lemma.\hfill$\diamondsuit$\\[0.2cm]
 \H Considering that (\ref{4z=0}) implies that $z$ vanishes in a neighborhood of the boundary of $B_{j,r^*_j}\cap\O$ (see the definition (\ref{4dr*j}) of $r^*_j$) , we can now construct the solution $\bar{z}_{j,\mu}$ of the problem $(\bar{P}_{j,\mu})$, which belongs to $H^1(\O)$, by
 \begin{equation}
 \bar{z}_{j,\mu}=\left\lbrace\begin{array}{ll}
 e^{-V(X(x,y),y)}(\int_{\a^{-1}(X(x,y))}^y \F{e^{V(X(x,y),t)}}{\scriptstyle\W U_2(X(x,y),t)}\bar{L}_{j,\mu}(X(x,y),t)\,dt)&\mr{if}\ (x,y)\in B_{j,r^*_j}\cap\O\\
 0& \mr{if}\ (x,y)\in\O\setminus B_{j,r^*_j},\end{array}\right.\end{equation}
 where the function $X$ is defined by (\ref{4chtv}), the functions $\bar{L}_{j,\mu}$, $U_2$ and $V$ are defined by (\ref{4dZLU2V}), the function $\a$ is defined by (\ref{4defaa-1}) and the real number $r^*_j$ by (\ref{4dr*j}), with a small enough real number $\mu$ verifying (\ref{4cmu}). Thus, thanks to (\ref{4dtpbp}), in this third case, we have proven that 
 \begin{equation}\label{zjmu3}
z_{j,\mu}\in H^1(\O)\ \mr{is\ the\ solution\ of\ the\ problem}\ (P_{j,\mu}).\end{equation}
 4) \tb{Fourth case}: $\mb{S^{-1}_j}\notin \gamma_{j-1}$, $\mb{S^{-1}_j}\in (E\cap S^c)$, $\mb{S^{1}_j}\in \gamma_{j+1}$.\\
\hspace*{2cm} \unitlength=1cm
\begin{picture}(8,5.5)
\put(4,4){\vector(1,0){2.5}}
\put(1,4){\vector(0,-1){2.5}}
\put(0.5,1.5){$x$}
\put(0.5,1){figure\ 3.7}
\put(6.1,3.7){$y$}
\color{red}
\qbezier(4,4)(2.5,4)(1,4)
\put(3.8,4.3){$\mb{S^1_j}$}
\put(0.8,4.3){$\mb{S^{-1}_j}$}
\qbezier(4,4)(5,3)(6,2)
\put(2,4.3){$\gamma_j$}
\put(5.6,2.9){$\gamma_{j+1}$}
\color{black}
\put(3,4.3){\vector(0,1){0.7}}
\put(5.3,3.2){\vector(1,1){0.5}}
\put(3.1,4.3){$\vec{n}_j$}
\put(4.7,3.5){$\vec{n}_{j+1}$}
\put(2.2,3.75){\vector(1,0){0.7}}
\put(2.4,3.3){$\vec{\tho}_-(\mb{S^{-1}_j})$}
\color{blue}
\qbezier(1,4)(0,4)(-1,4)
\qbezier(-1,4)(-1.5,3)(-2,2)
\put(-1,4.2){$\G^{0,+}$}
\color{green}
\put(0.8,3.8){\vector(-1,0){1}}
\put(-0.3,3.2){$\overrightarrow{u\scriptstyle ( S^{-1}_j)}$}
\end{picture}\\
\H The fourth case is not very different that the third case : $\mb{S^{-1}_j}$ is still the origin, the x-axis and the y-axis are defined in the same way, but, for $y<0$ small enough, the point $(0,y)$ belongs to $\Gamma^{0,+}$ and, therefore, we have $u_1(0,y)\le 0$. We denote by $\Gamma_{k_j}$ the side of the polygon which 
contains $\gamma_j$ and we set $\eta_j=d(\mb{S^{-1}_j},\P\Gamma_{k_j})$, where $d(.,.)$ is the euclidian distance in $\R^2$. According to the assumptions of the fourth case, $\eta_j>0$. Instead of (\ref{4dmu3}) which corresponds to the third case, we define $\mu'_3\le\min(\mu_2,\eta_j)$ such that
 \begin{eqnarray}
\forall\x\in B_{j,\mu_3'}\cap\O,\ \F{3}{2} u_2(0,0)\le u_2(\x)\le\F{1}{2} u_2(0,0)\ <0\ \ \mr{and}\nonumber\\ 
\forall y\in [-\mu'_3,\mu_3'],\ \di\F{1}{2}\di\F{\P u_1}{\P y}(0,0)y\le |u_1(0,y)|\le \di\F{3}{2}\di\F{\P u_1}{\P y}(0,0)y \label{4dmu3'}.\end{eqnarray}
Next, the proof of (\ref{4dI})-(\ref{4dmu4tB}) is slightly different in the case where $y=Y\le t\le 0$. First, we have 
$$X(0,t)=\int_0^t u_1(0,\theta)\,d\theta\le X(x,y).$$
Second, in view of (\ref{4dmu3'}), since $\int\limits_0^t u_1(0,\theta)\,d\theta\ge 0$, we still have
$$X(x,y)\le \F{3}{2}k(x,y)\ \mr{and}\ X(\F{\mu'_3}{2},t)\ge \F{\mu'_3}{4}|u_2(0,0)|$$
and (\ref{4dZLU2V}), (\ref{4dI}) and (\ref{4dmu4tB}) are still verified with $\mu'_3$ in the place of $\mu_3$.
In the same way, (\ref{4defaa-1}) and (\ref{4cXY}) run unchanged, while (\ref{4zsol}) is verified with $K'$ in the place of $K$ where $K'$ is defined by
\begin{equation}\label{4dK'}
K'=\min(\F{\mu_3'}{6},\mu_4,\mu_5,|\gamma_j|)
.\end{equation}
 Afterwards, the case where $0\le y\le t\le \a^{-1}(X(x,y))$ remains unchanged and we still have (\ref{4z=0}) for $y\ge 0$, with $r'_{1,j}$ and $r'_{2,j}$ in the place of $r_{1,j}$ and $r_{2,j}$, where $r'_{1,j}$ and $r'_{2,j}$ are defined by
 \begin{equation}\label{4r'12}
  r'_{1,j}=\min(\F{|u_2(0,0)|K'}{12},\F{|\F{\P u_1}{\P y}(0,0)| K'^2}{288})\ \mr{and}\ r'_{2,j}=\F{K'|u_2(0,0)|}{6}.\end{equation}
 \H When $y\le t\le\a^{-1}(X(x,y))$, with $y<0$, we consider first $0\le t\le \a^{-1}(X(x,y))$ and second $y\le t<0$.\\
 \H If $0\le t\le \a^{-1}(X(x,y))$ and $(x,y)\in B_{j,\F{K'}{6}}^+$, we have $t\le \F{\mu'_3}{6}$ and
$$\a(t)=X(0,t)\le X(x,y)\le \F{3}{2}k(x,y).$$ Applying (\ref{4ikB}) with $r=|u_2(0,0)|\F{K'}{6}$, we obtain $k(x,y)\le |u_2(0,0)|\di\F{K'}{6}$.
Hence, we derive
$$X(0,t)\le X(x,y)\le |u_2(0,0)|\di\F{K'}{4}\le X(\F{K'}{2},t).$$
Then, there exists $x_t\in [0,\F{K'}{2}]$ such that
$$X(x,y)=X(x_t,t)\ \mr{with}\ (x_t,t)\in B_{j,\mu'_3}^+.$$
Next, as previously, we derive that $\forall(x,y)\in C(\mb{S^{-1}_j},2\sqrt{\F{r'_{1,j}}{|\F{\P u_1}{\P y}(0,0)|}},\F{r'_{2,j}}{|u_2(0,0)|})\cap\ov{\O}$ with $y< 0,$
$$\int_{\a^{-1}(X(x,y))}^0 \F{e^{V(X(x,y),t)}}{\scriptstyle\W U_2(X(x,y),t)}\bar{L}_{j,\mu}(X(x,y),t)\,dt=0.$$
Finally, for the case where $y\le t\le 0$, we process in the same way as previously and we obtain that (\ref{4z=0}) is verified for $y<0$ with $r'_{1,j}$ and $r'_{2,j}$ in the place of $r_{1,j}$ and $r_{2,j}$.\\
\H The rest of the proof is the same as in the third case. Thus, in the fourth case, we have proven that 
 \begin{equation}\label{zjmu4}
z_{j,\mu}\in H^1(\O)\ \mr{is\ the\ solution\ of\ the\ problem}\ (P_{j,\mu}).\end{equation}
\H By localization, all the other cases, where $\mb{S^{-1}_j}\in E$ and (or) $\mb{S^{1}_j}\in E$, can be solved as in the third case or the fourth case.
\hfill$\diamondsuit$\\
\section{Appendix}
 \H In the two following examples, the domains $\O$ are no longer a bounded polygon, but domains of class $C^{1,1}$. Even if in this article, we mainly deal with bounded polygons, it seems to us interesting to show that the regularity of the solution $z$ of the transport problem in domains of class $C^{1,1}$ seems still linked to the multiplicity of the roots of the equation $\u\,.\,\mb{n}=0$ at the end-points of $\G^-$, in the case where $\u\,.\,\tho_-$ is negative.
\subsection{Example 6 : $\O=C(I(0,1),0.5)$,\ $l(x,y)=1$, $\u(x,y)=(x,-y)$.}
\begin{minipage}[b]{11cm}
\H In this example, the boundary, which is the circle of center $I(0,1)$ and of radius $R=0.5$, is very regular, but the function $\u.\mb{n}$ vanishes at the boundary points of $\G^-$, which leads to a discontinuity for the 
partial derivatives of the solution $z$ in these points. \\
\H The equation of $\G$ is $\left\lbrace\begin{array}{l}x=\di\F{1}{2}\cos t\\
y=1+\di\F{1}{2}\sin t\end{array}\right.,\ t\in ]-\pi,\pi]$ and the unit exterior normal is $\mb{n}=(\cos t,\sin t)$. Let us determine the sets $\G^-$, $\G^0$ and $\G^+$.
On $\G$, we have $$(\u\,.\,\mb{n})(t)=-\sin^2 t-\sin t+\F{1}{2},$$ that vanishes for $t_0=\arcsin(\F{\sqrt{3}-1}{2})$ and $t_1=\pi-\arcsin(\F{\sqrt{3}-1}{2})$, and 
$\G^-$ is the open arc of the circle $\G$ defined by $t_0<t<t_1$, $\G^+$ is $\G\setminus\ov{\G^-}$, $\G^0=\emptyset$.
\end{minipage}
\hspace*{1cm}
\unitlength=1cm
\begin{picture}(6,6)
\put(1.1,1.5){\vector(0,1){4.5}}\put(-1,2){\vector(1,0){5}}
\linethickness{0.03cm}
\color{red}
\qbezier(1.93,4.365)(1.8776,4.4794)(1.8253,4.5646)
\qbezier(1.8253,4.5646)(1.8196,4.5729)(1.8160,4.577)
\qbezier(1.8160,4.577)(1.6967,4.7174)(1.6216,4.7833)
\qbezier(1.6216,4.7833)(1.4536,4.8912)(1.2675,4.9636)
\qbezier(1.2675,4.9636)(1.1502,4.9887 )(1,5)
\qbezier(1,5)(0.89101,4.994)(0.7728,4.9738)
\qbezier(0.7728,4.9738)(0.67671,4.9463)(0.58385,4.9093)
\qbezier(0.58385,4.9093)(0.49515,4.8632)(0.4115,4.8085)
\qbezier(0.4115,4.8085)(0.29729,4.7115)(0.19886,4.5985)
\qbezier(0.19886,4.5985)(0.14311,4.5155)(0.11374,4.4632)
\qbezier(0.11374,4.4632)(0.08756,4.4092)(0.069,4.367)
\color{blue}
\qbezier(0.069,4.367)(0.04221,4.2875)(0.01001,4.1411)
\qbezier(0.01001,4.1411)(0,4)(0.00171,3.9416)
\qbezier(0.00171,3.9416)(0.02164,3.7931)(0.06354,3.6492)
\qbezier(0.06354,3.6492)(0.12648,3.5132)(0.20903,3.3881)
\qbezier(0.20903,3.3881)(0.30935,3.2768)(0.42518,3.1817)
\qbezier(0.42518,3.1817)(0.55391,3.105)(0.69267,3.0484)
\qbezier(0.69267,3.0484)(0.83832,3.0132)(0.98761,3.0001)
\qbezier(0.98761,3.0001)(1.1372,3.0095)(1.2837,3.0411)
\qbezier(1.2837,3.0411)(1.4238,3.0942)(1.5544,3.1677)
\qbezier(1.5544,3.1677)(1.6725,3.2599)(1.7756,3.3687)
\qbezier(1.7756,3.3687)(1.8612,3.4917)(1.9275,3.6261)
\qbezier(1.9275,3.6261)(1.9729,3.7689)(1.9965,3.9169)
\qbezier(1.9965,3.9169)(2,4)(1.9978,4.0668)
\qbezier(1.9978,4.0668)(1.9861,4.166)(1.9646,4.2637)
\qbezier(1.9646,4.2637)(1.9471,4.321)(1.9305,4.366)
\color{black}
\put(0.4,5.05){1.5}
\put(0.6,4){1}
\put(1.1,3.85){I}
\put(0.4,2.7){0.5}
\put(0.9,4){\line(1,0){0.1}}
\put(3,1.95){\line(0,1){0.1}}
\put(2.9,1.6){1}
\put(0.6,6){y}
\put(4,1.7){x}
\put(0.6,1.65){O}
\thinlines
\qbezier[15](2,2)(2,3)(2,4)
\qbezier[15](0,2)(0,3)(0,4)
\put(1.7,1.6){0.5}
\put(-0.4,1.6){-0.5}
\color{red}
\put(1.5,5){$\G^-$}
\color{blue}
\put(1.5,2.8){$\G^+$}
\color{black}
\put(1,0){figure\ 4.8}
\end{picture}\\[0.2cm]
\H Note that, we can easily verify that $(\u\,.\,\tho_{-})(t)$ is negative for $t=t_0$ and $t=t_1$, so that assumptions analogous to the assumptions of (\ref{4ctheun})
are verified  at points where $\u\,.\,\mb{n}$ va-\linebreak nishes in this example. As in the first three examples, we have 
\begin{equation}\label{2ex4z}
\forall (x,y)\in\O,\ z(x,y)=1-2y+y\,C(xy),\ 
 \forall(x,y)\in\G^-,\ C(xy)=2-\di\F{1}{y}.\end{equation} 
 Setting $X=xy$, we must compute the function $\alpha$ such that $y=\alpha(X)$, for $(x,y)\in \G^-$.\\
 1) \tb{First case} : $t_0< t\le \F{\pi}{2}$, $(x,y)\in\G^-\cap \R_+^2$. \\
 $x=\F{1}{2}\sqrt{1-(2y-2)^2}=\F{1}{2}\sqrt{(2y-1)(3-2y)}$, which imply $X=\di\F{y}{2}\sqrt{(2y-1)(3-2y)}$. We compute $y(t_0)=\F{3+\sqrt{3}}{4}$ and $X(t_0)=\F{\sqrt{2\sqrt{3}}(3+\sqrt{3})}{16}$. Considering the function 
 \begin{equation}\label{2dg} g:\ y\mapsto X=\F{y}{2}\sqrt{(2y-1)(3-2y)},\end{equation} which is defined on the set $[\F{1}{2},\F{3}{2}]$. Since we have, $\forall y\in ]\F{1}{2},\F{3}{2}[,\ g'(y)=-\di\F{4(y-\F{3+\sqrt{3}}{4})(y-\F{3-\sqrt{3}}{4})}{\sqrt{(2y-1)(3-2y)}}$,\\[0.2cm] the statement of changes of $g$ is 
  \hspace*{0.2cm} \begin{tabular}{|c|l|}
\hline
$y$ &$\F{1}{2}$\hspace*{2.2cm}$\F{3+\sqrt{3}}{4}$\hspace*{2.3cm} $\F{3}{2}$\\
\hline
$g'(y)$&$\|$\,$+\iy$\hspace*{0.7cm}+\hspace*{0.7cm}0\hspace*{0.7cm}\scalebox{2}[1]{-}\hspace*{1cm}$-\iy$\,$\|$\\
\hline
$g$ &\begin{picture}(4,1.5)
\put(0,0.2){0}
\put(0.3,0.6){\vector(4,1){1.7}}
\put(2.5,1){$X(t_0)$}
\put(3.7,1){\vector(4,-1){1.7}}
\put(5.73,0.3){0}
\end{picture}
\\ \hline
\end{tabular}\ .\\[0.2cm]
Since $g$ is strictly decreasing from $[y(t_0),\F{3}{2}]$ to 
 $[0,X(t_0)]$, therefore $g_{|[y(t_0),\F{3}{2}]}$ has an inverse function and we have
  \begin{equation}\label{2dalpha1}
  \alpha_{|[0,X(t_0)]}=(g_{|[y(t_0),\F{3}{2}]})^{-1}. 
  \end{equation}
 Finally, we obtain, $\mr{for}\ (x,y)\in \G^-\cap \R_+^2,\ y=\alpha(X)\Longleftrightarrow y =g(X)$.\\
 2) \tb{Second case} : $ \F{\pi}{2}\le t\le \pi-t_0$, $(x,y)\in\G^-\cap \R_-\t\R_+$. \\
In the same way, we have $X=-\di\F{y}{2}\sqrt{(2y-1)(3-2y)}$ and we define $\alpha$ on $[-X(t_0),0]$ by
\begin{equation}\label{2dalpha2}
  \alpha_{|[-X(t_0),0]}=((-g)_{|[y(t_0),\F{3}{2}]})^{-1}. 
  \end{equation}
  The statement of changes of the even function $\alpha$ is :
  \begin{tabular}{|c|l|}
\hline
$X$ &\hspace*{-0.2cm}$-X(t_0)$\hspace*{1.6cm}0\hspace*{1.9cm} $+X(t_0)$\\
\hline
$\alpha'(X)$&$\|$\,$+\iy$\hspace*{0.7cm}+\hspace*{0.7cm}0\hspace*{0.7cm}\scalebox{2}[1]{-}\hspace*{1cm}$-\iy$\,$\|$\\
\hline
$\alpha$ &\begin{picture}(4,1.5)
\put(-0.1,0.2){$\F{3+\sqrt{3}}{4}$}
\put(0.8,0.6){\vector(4,1){1.7}}
\put(2.75,1){$\F{3}{2}$}
\put(3.2,1){\vector(4,-1){1.7}}
\put(5,0.3){$\F{3+\sqrt{3}}{4}$}
\end{picture}
\\ \hline
\end{tabular}\\[0.2cm]
From (\ref{2ex4z}), we derive the solution of the example 6
\begin{equation}\label{2zex4}
\forall (x,y)\in \O,\ z(x,y)= 1-\F{y}{\a(xy)}.\end{equation}
Let us show that $z$ belongs to $H^1(\O)$. We compute 
$$\forall (x,y)\in \O\setminus\{(\pm x(t_0),y(t_0))\},
\ z'_x(x,y)=\F{y^2\a'(xy)}{(\a(xy))^2},\ z'_y(x,y)=\F{xy\a'(xy)}{(\a(xy))^2}-\F{1}{\a(xy)},$$
with $x(t_0)=\di\F{\sqrt{2\sqrt{3}}}{4}$ and $y(t_0)=\F{3+\sqrt{3}}{4}$. Hence, we derive
\begin{equation}\label{2a'zH1}
z\ \mr{belongs\ to}\ H^1(\O)\Longleftrightarrow \int\int_{\O}(\a'(xy))^2\,dxdy< +\iy.\end{equation}

\H Let us set $\O_+=\O\cap\R_+^2$ and note that $\int\int_{\O}(\a'(xy))^2\,dxdy=2\int\int_{\O_+}(\a'(xy))^2\,dxdy.$ In order to show that the last integral converges, we split $\O_+$ in two subdomains: $$\O_+^1=\O_+\cap (\R_+\t [\F{3+\sqrt{3}}{4},\F{3}{2}])\ \mr{and}\ \O_+^2=\O_+\cap (\R_+\t [\F{1}{2},\F{3+\sqrt{3}}{4}]).$$
 1) We compute $\int\int_{\O_+^1}(\a'(xy))^2\,dxdy$ by making the substitution $\left\lbrace\begin{array}{l}u=\a(xy)\\
y=y\end{array}\right.\ \Longleftrightarrow \left\lbrace\begin{array}{l}x=\F{g(u)}{y}\\y=y\end{array}\right.$, the jacobian of which is $\di\F{g'(u)}{y}$. Since 
$\forall u\in ]\F{3+\sqrt{3}}{4},\F{3}{2}],\ \a'(g(u))=\di\F{1}{g'(u)}$, we obtain
$$\int\int_{\O_+^1}(\a'(xy))^2\,dxdy=\int_{\F{3+\sqrt{3}}{4 }}^{\F{3}{2}}\F{1}{y}\,dy(\int_y^{\F{3}{2}}\F{1}{|g'(u)|}\,du).$$
Considering that, $\forall u\in ]\F{3+\sqrt{3}}{4},1.5[,\ g'(u)=-\di\F{4(u-\F{3+\sqrt{3}}{4})(u-\F{3-\sqrt{3}}{4})}{\sqrt{(2u-1)(3-2u)}}$, 
we can verify
\begin{equation}\label{2bO1}\int\int_{\O_+^1}(\a'(xy))^2\,dxdy\le K_1\int_{\F{3+\sqrt{3}}{4 }}^{\F{3}{2}}dy(\int_y^{\F{3}{2}}\F{1}{u-\F{3+\sqrt{3}}{4}}\,du)=\F{3-\sqrt{3}}{4}K_1<+\iy,\end{equation} 
with $K_1=\di\F{\sqrt{3}-1}{3}$. \\
2) For $\int\int_{\O_+^2}(\a'(xy))^2\,dxdy$, it is more complicated. Setting $X=xy$, we obtain
$$\int\int_{\O_+^2}(\a'(xy))^2\,dxdy=\int_{\F{1}{2}}^{\F{3+\sqrt{3}}{4}}\,dy(\int_0^{\F{g(y)}{y}}(\a'(xy))^2\,dx)=
\int_{\F{1}{2}}^{\F{3+\sqrt{3}}{4}}\F{1}{y}\,dy(\int_0^{g(y)}(\a'(X))^2\,dX),$$
where the function $g$ is defined by (\ref{2dg}). Next, making the substitution $X=g(u)$, for $u\in [\F{1}{2},\F{3+\sqrt{3}}{4}]$, we derive
$$\int\int_{\O_+^2}(\a'(xy))^2\,dxdy=\int_{\F{1}{2}}^{\F{3+\sqrt{3}}{4}}\F{1}{y}\,dy(\int_{\F{1}{2}}^y(\a'(g(u)))^2g'(u)\,du).$$
However, the complication comes from the fact that, for $\F{1}{2}\le u<\F{3+\sqrt{3}}{4}$,\ $\a(g(u))\not= u$ since $\a(g(u))> \F{3+\sqrt{3}}{4}$.
Let us define the function $\beta$ on the set $[ \F{1}{2},\F{3+\sqrt{3}}{4}]$ by
$$\forall u\in [\F{1}{2},\F{3+\sqrt{3}}{4}],\ \beta(u)=\a(g(u)).$$
Since $g(\beta(u))=g(u)$, then, for $\F{1}{2}<u<\F{3+\sqrt{3}}{4}$, 
$$g'(u)=g'(\beta(u))\beta'(u)\ \mr{and}\ \a'(g(u))=\a'(g(\beta(u)))=\F{1}{g'(\beta(u))}=\F{\beta'(u)}{g'(u)}.$$
Hence, we can write
$$\int\int_{\O_+^2}(\a'(xy))^2\,dxdy=\int_{\F{1}{2}}^{\F{3+\sqrt{3}}{4}}\F{1}{y}\,dy(\int_{\F{1}{2}}^y\F{(\beta'(u))^2}{g'(u)}\,du).$$
Let us show that $\beta'$ is bounded on the set $]\F{1}{2},\F{3+\sqrt{3}}{4}[$ and that we can extend $\beta'$ on $[\F{1}{2},\F{3+\sqrt{3}}{4}]$ by continuity.
Computing $\beta'(u)$ yields
$$\beta'(u)=\F{u(u-\F{3+\sqrt{3}}{4})(u-\F{3-\sqrt{3}}{4})}{\beta(u)(\beta(u)-\F{3+\sqrt{3}}{4})(\beta(u)-\F{3-\sqrt{3}}{4})}.$$
Note that the right hand previous expression extends $\beta'$ by continuity in $\F{1}{2}$. It remains to compute the limit of $\beta'$ in $\F{3+\sqrt{3}}{4}$.
Applying to the function $g$ the Taylor-Lagrange formula in the neighborhood of $\F{3+\sqrt{3}}{4}$, we obtain
$$g(u)=g(\F{3+\sqrt{3}}{4})+\F{1}{2}g''(c)(u-\F{3+\sqrt{3}}{4})^2=g(\beta(u))= g(\F{3+\sqrt{3}}{4})+\F{1}{2}g''(d)(\beta(u)-\F{3+\sqrt{3}}{4})^2,$$
with $c\in [u,\F{3+\sqrt{3}}{4}]$ and $d\in [\F{3+\sqrt{3}}{4},\beta(u)]$, which implies
$\lim\limits_{u\to \F{3+\sqrt{3}}{4}^-}\F{(u-\F{3+\sqrt{3}}{4})^2}{(\beta(u)-\F{3+\sqrt{3}}{4})^2}=1$ and 
$\lim\limits_{u\to \F{3+\sqrt{3}}{4}^-}\F{u-\F{3+\sqrt{3}}{4}}{\beta(u)-\F{3+\sqrt{3}}{4}}=-1$, since $\F{u-\F{3+\sqrt{3}}{4}}{\beta(u)-\F{3+\sqrt{3}}{4}}\le 0$. Hence, we derive $\lim\limits_{u\to \F{3+\sqrt{3}}{4}^-}\beta'(u)=-1$ and, therefore, there exists a constant $K_2>0$, such that, for $u\in [\F{1}{2},\F{3+\sqrt{3}}{4}]$, $|\beta'(u)|\le K_2$.
Then, we have $$\int\int_{\O_+^2}(\a'(xy))^2\,dxdy\le (\sqrt{3}+1)K_2^2\int_{\F{1}{2}}^{\F{3+\sqrt{3}}{4}}\,dy(\int_{\F{1}{2}}^y\F{1}{\F{3+\sqrt{3}}{4}-u}\,du)\le
(\F{(\sqrt{3}+1)K_2}{2})^2<+\iy$$ and, with (\ref{2a'zH1} and (\ref{2bO1}), we derive that the solution $z$ belongs to 
$H^1(\O)$. Finally, although $\u\,.\,\mb{n}$ vanishes at the end points of $\G^-$, the problem (\ref{pbh2}) is well-posed, probably because the function $\u\,.\,\mb{n}$ has 
only simple roots at the end points of $\G^-$ with, in addition, $\u\,.\,\tho_-$ negative at these end points.
\subsection{\tb{Example 7 : $\O=\O_7$, $l(x,y)=1$, $\u(x,y)=(x,-y)$.}}
\begin{minipage}[b]{10cm}
\H The boundary of $\O_7$ is composed of two half semicircles, linked up by two segments (see the figure 4.9). The boundary is of class $C^{1,1}$
 but the arc of circle $\G^-$ is adjacent to the segment $\G^0$, which leads to a discontinuity for the 
partial derivatives of the solution $z$. But, as in the example 4, \tb{this discontinuity is such that the solution} $\mb{z}$ \tb{does not belong to} $\mb{H^1(\O)}$,
\tb{as we shall see further}. The parametric equation of the upper semicircle is 
$\left\lbrace\begin{array}{l}x=\di\F{1}{2}+\di\F{1}{2}\cos t\\[0.2cm]
y=1+\di\F{1}{2}\sin t\end{array}\right.,\ t\in [0,\pi]$. 
\end{minipage}\hspace*{3.5cm}
\unitlength=1cm
\begin{picture}(6,3)
\put(0.13,1.5){\vector(0,1){4.3}}
\put(-1,2){\vector(1,0){3.7}}
\linethickness{0.03cm}
\color{blue}
\qbezier(2,4)(1.998,4.02)(1.9992,4.04)
\qbezier(1.9992,4.04)(1.995,4.0998)(1.9888,4.1494)
\qbezier(1.9888,4.1494)(1.9689,4.2474)(1.9492,4.3146)
\qbezier(1.9492,4.3146)(1.9394,4.3429)(1.93,4.365)
\qbezier(1.93,4.365)(1.8776,4.4794)(1.8253,4.5646)
\qbezier(1.8253,4.5646)(1.8196,4.5729)(1.8170,4.577)
\color{red}
\qbezier(1.8170,4.577)(1.6967,4.7174)(1.6216,4.7833)
\qbezier(1.6216,4.7833)(1.4536,4.8912)(1.2675,4.9636)
\qbezier(1.2675,4.9636)(1.1502,4.9887 )(1,5)
\qbezier(1,5)(0.89101,4.994)(0.7728,4.9738)
\qbezier(0.7728,4.9738)(0.67671,4.9463)(0.58385,4.9093)
\qbezier(0.58385,4.9093)(0.49515,4.8632)(0.4115,4.8085)
\qbezier(0.4115,4.8085)(0.29729,4.7115)(0.19886,4.5985)
\qbezier(0.19886,4.5985)(0.14311,4.5155)(0.11374,4.4632)
\qbezier(0.11374,4.4632)(0.08756,4.4092)(0.069,4.367)
\qbezier(0.069,4.367)(0.04221,4.2875)(0.01001,4.1411)
\qbezier(0.01001,4.1411)(0,4)(0.00171,3.9416)
\color{green}
\linethickness{0.03cm}
\qbezier(0,4)(0,3.75)(0,3.5)
\put(-0.6,3.6){$\G^0$}
\color{blue}
\qbezier(0.00171,3.4416)(0.02164,3.2931)(0.06354,3.1492)
\qbezier(0.06354,3.1492)(0.12648,3.0132)(0.20903,2.8881)
\qbezier(0.20903,2.8881)(0.30935,2.7768)(0.42518,2.6817)
\qbezier(0.42518,2.6817)(0.55391,2.605)(0.69267,2.5484)
\qbezier(0.69267,2.5484)(0.83832,2.5132)(0.98761,2.5001)
\qbezier(0.98761,2.5001)(1.1372,2.5095)(1.2837,2.5411)
\qbezier(1.2837,2.5411)(1.4238,2.5942)(1.5544,2.6677)
\qbezier(1.5544,2.6677)(1.6725,2.7599)(1.7756,2.8687)
\qbezier(1.7756,2.8687)(1.8612,2.9917)(1.9275,3.1261)
\qbezier(1.9275,3.1261)(1.9729,3.2689)(1.98,3.35)
\qbezier(1.98,3.35)(1.9965,3.4169)(2,3.5)
\qbezier(2,3.5)(2,3.75)(2,4)
\color{black}
\put(-0.6,5.05){1.5}
\put(-0.4,4){1}
\put(1.1,3.85){$\O_7$}
\put(-0.8,3.2){0.75}
\put(-0.8,2.4){0.25}
\qbezier[10](0,2.5)(0.5,2.5)(1,2.5)
\qbezier[10](0,5)(0.5,5)(1,5)
\qbezier[7](1,2)(1,2.25)(1,2.5)
\qbezier[15](2,2)(2,2.5)(2,3.5)
\put(1.9,1.6){1}
\put(-0.4,5.7){y}
\put(2.3,1.7){x}
\put(-0.4,1.65){O}
\thinlines
\put(0.7,1.6){0.5}
\color{red}
\put(1,5.2){$\G^-$}
\color{blue}
\put(2.1,2.7){$\G^+$}
\color{black}
\put(0,0){figure\ 4.9}
\end{picture}\\
\H Let us determine the sets $\G^-$, $\G^0$ and $\G^+$. Again, we have $\mb{n}=(\cos t,\sin t)$.
On the upper semicircle, we compute $(\u\,.\,\mb{n})(t)=\cos (\F{t}{2})(\cos (\F{t}{2})-2\sin(\F{t}{2})(\sin t+1))$. In view of 
$\cos (\F{t}{2})>0$ for $0\le t<\pi$, $(\u\,.\,\mb{n})(t)$ has the same sign that $f(t)=\cos (\F{t}{2})-2\sin(\F{t}{2})(\sin t+1).$\\
1) For $\F{\pi}{2}\le t\le \pi$, $f(t)=\cos (\F{t}{2})(1-4\sin^2(\F{t}{2}))-2\sin(\F{t}{2})<0$.\\
2) For $0\le t\le \F{\pi}{2}$, $f(t)=\cos^3 (\F{t}{2})(-2\tan^3(\F{t}{2})-3\tan^2(\F{t}{2})-2\tan(\F{t}{2})+1)$. Setting $\theta=\tan(\F{t}{2})$, we must 
study the sign of the polynomial $g(\theta)= -2\theta^3 -3\theta^2-2\theta+1$, for $0\le \theta\le 1$. We have $g'(\theta)=-6\theta^2-6\theta-2<0$, 
$g(0)=1$, $g(1)=-6$. Then the continuity and the strict decreasing of $g$ implies that there exists an unique number $\theta_0\in ]0,1[$, such that $g(\theta_0)=0$.
Finally, for $0\le t\le \pi$, $(\u\,.\,\mb{n})(t)$ vanishes for two values : 
\begin{equation}\label{2dt0}t_0=2\arctan(\theta_0)\approx 0.614\ \mr{and}\ \pi\end{equation} and 
$$\G^-\ \mr{is\ the\ open\ arc\ of\ circle\ defined\ by}\ \left\lbrace\begin{array}{l}x=\di\F{1}{2}+\di\F{1}{2}\cos t\\[0.2cm]
y=1+\di\F{1}{2}\sin t\end{array}\right.,\ t\in ]t_0,\pi[,$$ while the part of the previous semicircle, defined by $t\in [0,t_0[$, is included into $\G^+$. 
Next, the part of $\G$, defined by $\left\lbrace\begin{array}{l}x=0\\
\F{3}{4}< y< 1\end{array}\right.$ is $\G^0$ and we are going to show that the lower semicircle is included in $\G^+$. Indeed, the parametric equations of the 
lower semicircle is $\left\lbrace\begin{array}{l}x=\di\F{1}{2}+\di\F{1}{2}\cos t\\[0.2cm]
y=\di\F{3}{4}+\di\F{1}{2}\sin t\end{array}\right.,\ t\in [-\pi,0]$. Considering that, for $t\in [-\pi,0]$,
$$(\u\,.\mb{n})(t)=\F{1}{2}\cos t +\F{1}{2}\cos^2 t-\F{3}{4}\sin t- \F{1}{2}\sin^2 t,$$
we distinguish two cases :\\
1) $t\in ]-\pi,-\F{3\pi}{4}[\cup]-\F{\pi}{2},0]$\\
$(\u\,.\mb{n})(t)=\F{1}{2}(\cos t-\sin t)(1+\cos t+\sin t)-\F{1}{4}\sin t >0$.\\
2) $t\in [-\F{3\pi}{4},-\F{\pi}{2}]$\\[0.1cm]
$\F{1}{2}\cos t-\F{1}{2}\sin^2 t=-\F{5}{8}+\F{1}{2}(\cos t+\F{1}{2})^2\ge -\F{5}{8}$ and $\F{1}{2}\cos^2 t-\F{3}{4}\sin t=\F{25}{32}-\F{1}{2}(\sin t+\F{3}{4})^2\ge 
\F{3}{4}$. Therefore, $(\u\,.\mb{n})(t)\ge \F{1}{8}>0$.\\
Finally, $\left\lbrace\begin{array}{l}x=1\\
\F{3}{4}< y< 1\end{array}\right.$ is included in $\G^+$, which ends the determining of the sets $\G^-$, $\G^0$ and $\G^+$, see figure 4.9.\\
\H Again, it is easy to verify that $(\u\,.\,\tho_-)(t)$ is negative for $t=t_0$ and $t=\pi$. As in example 6, we have 
\begin{equation}\label{2ex5z}
\forall (x,y)\in\O,\ z(x,y)=1-2y+y\,C(xy),\ 
 \forall(x,y)\in\G^-,\ C(xy)=2-\di\F{1}{y}.\end{equation} 
 Setting $X=xy$, we must compute the function $\alpha$ such that $y=\alpha(X)$, for $(x,y)\in \G^-$.\\
 1) \tb{First case} : $(x,y)\in\G^-\cap ([0,\F{1}{2}]\t\R_+)$. \\
 $X=\di\F{y}{2}(1-\sqrt{(2y-1)(3-2y)}\,)$. Considering the function
\begin{equation}\label{2dg5} g:\ y\mapsto X=\F{y}{2}(1-\sqrt{(2y-1)(3-2y)}).\end{equation}
We have, for $\F{1}{2}<y< \F{3}{2}$,
\begin{equation}\label{2g'5}
g'(y)=\F{8y^2-12y+3+\sqrt{(2y-1)(3-2y)}}{2\sqrt{(2y-1)(3-2y)}}.\end{equation}
For $\di\F{1}{2}<y< 1$, $8y^2-12y+3< -1$, $\sqrt{(2y-1)(3-2y)}<1$, therefore, $g'(y)<0$. Moreover, $g'(1)=0$.\\
For $\di\F{3+\sqrt{3}}{4}<y<\F{3}{2}$, $8y^2-12y+3>0$, therefore $g'(y)>0$.\\
For $1<y<\F{3+\sqrt{3}}{4}$, we note that $g'(y)$ has the same sign that 
$$h(y)=8y^2-12y+3+\sqrt{(2y-1)(3-2y)}$$ 
and $h'(y)=16y-12-\di\F{4(y-1)}{\sqrt{(2y-1)(3-2y)}}$. In view of $\di\F{4(y-1)}{\sqrt{(2y-1)(3-2y)}}\le\F{\sqrt{6}-\sqrt{2}}{3^{0.25}}<1$ and $(16y-12)\ge 4$,  we derive $h'(y)>0$. Since $h(1)=0$, we obtain that $h(y)>0$, that is to say, $g'(y)>0$.\\
 Thus, we obtain the statement of changes of $g$ :
\hspace*{0.2cm} \begin{tabular}{|c|l|}
\hline&\\[-0.35cm]
$y$ &$\F{1}{2}$\hspace*{2.48cm}$1$\hspace*{2.65cm} $\F{3}{2}$\\[0.1cm]
\hline
$g'(y)$&$\|$\,$-\iy$\hspace*{0.7cm}\scalebox{2}[1]{-}\hspace*{0.7cm}0\hspace*{0.7cm}+\hspace*{1cm}$+\iy$\,$\|$\\
\hline
$g$ &\begin{picture}(4,1)
\put(0,0.6){$\F{1}{4}$}
\put(0.5,0.6){\vector(4,-1){1.7}}
\put(2.7,0){0}
\put(3.5,0.2){\vector(4,1){1.7}}
\put(5.73,0.6){$\F{3}{4}$}
\end{picture}
\\ \hline
\end{tabular}\\[0.2cm]
Finally, the function $g$, defined by (\ref{2dg5}), is strictly increasing from $[1,\F{3}{2}]$ to 
 $[0,\F{3}{4}]$, therefore $g_{|[1,\F{3}{2}]}$ has an inverse function and we define $\a$ on the set $[0,\F{3}{4}]$ by
  \begin{equation}\label{2dalpha5}
  \alpha_{|[0,\F{3}{4}]}=(g_{|[1,\F{3}{2}]})^{-1}, 
  \end{equation} that verifies $y=\alpha(X)\Longleftrightarrow X =g(y)\ \mr{for}\ (x,y)\in \G^-\cap ([0,\F{1}{2}]\t \R_+).$ \\
2) \tb{Second case} : $(x,y)\in\G^-\cap ([\F{1}{2},1]\t\R_+)$. \\
 $X=\di\F{y}{2}(1+\sqrt{(2y-1)(3-2y)}\,)$, with $y(t_0)<y<\F{3}{2}$, where $t_0$ is defined by (\ref{2dt0}) and $y(t_0)\approx 1.288$. Considering the function
$$ \tilde{g}:\ y\mapsto X=\F{y}{2}(1+\sqrt{(2y-1)(3-2y)}).$$ We compute $\tilde{g}'(y)=\di\F{-8y^2+12y-3+\sqrt{(2y-1)(3-2y)}}{2\sqrt{(2y-1)(3-2y)}}$. We can verify that
$$\forall t\in [0,\F{\pi}{2}],\ (\u\,.\,\mb{n})(t)=0\Longleftrightarrow -8y^2+12y-3+\sqrt{(2y-1)(3-2y)}=0.$$
Hence, we derive that $\tilde{g}'$ vanishes at $y(t_0)$ and, since the numerator of $\tilde{g}'$ is a strictly decrea\linebreak sing function on $[y(t_0),\F{3}{2}]$, we obtain that 
$\tilde{g}'$ is strictly negative on $]y(t_0),\F{3}{2}[$, which implies that the function $\tilde{g}$ is strictly decreasing 
from $[y(t_0), \F{3}{2}]$ to $[\F{3}{4},\tilde{g}(y(t_0))]$, with $\tilde{g}(y(t_0))\approx 1.17$. Therefore $\tilde{g}_{|[y(t_0), \F{3}{2}]}$ has an inverse function and we define $\a$ on the set $[\F{3}{4},\tilde{g}(y(t_0))]$ by
  \begin{equation}\label{2dalpha5b}
  \alpha_{|[\F{3}{4},\tilde{g}(y(t_0))]}=(\tilde{g}_{|[y(t_0), \F{3}{2}]})^{-1}, 
  \end{equation} that verifies $y=\alpha(X)\Longleftrightarrow X =\tilde{g}(y)\ \mr{for}\ (x,y)\in \G^-\cap ([\F{1}{2},1]\t \R_+).$ \\
  Finally, from (\ref{2ex5z}), 
  we derive the solution of example 7 :
\begin{equation}\label{2zex5}
\forall (x,y)\in \O,\ z(x,y)= 1-\F{y}{\a(xy)},\end{equation} where the function $\a$ is defined on the interval $[0,\tilde{g}(y(t_0))]$ by (\ref{2dalpha5}) and (\ref{2dalpha5b}).\\
\H Let us show that $z$ does not belongs to $H^1(\O)$. We again compute
$$\forall (x,y)\in \O\setminus\{(0,1))\},
\ z'_x(x,y)=\F{y^2\a'(xy)}{(\a(xy))^2}.$$
We integrate $(z'_x(x,y))^2$ on a domain $\O^*=\{(x,y)\in\R^2,\ \di\F{3}{4}\le y\le 1,\ 0\le x\le \di\F{1}{4y}\}$, which is included in $\O$. Moreover, since 
in $\O^*$, we have $0\le xy\le \di\F{1}{4}\le \di\F{3}{4}$, we use the expression of the fonction $\a$ defined by (\ref{2dalpha5}) 
from the function $g$, defined by (\ref{2dg5}). Clearly, $$\int\int_{\O}(z'_x(x,y))^2\,dxdy\ge\int\int_{\O^*}(z'_x(x,y))^2\,dxdy.$$ Let us show that the integral $\int\int_{\O^*}(z'_x(x,y))^2\,dxdy =+\iy$. 
First, we compute $\int\int_{\O_*}(z'_x(x,y))^2\,dxdy$ by making the substitution $\left\lbrace\begin{array}{l}X=xy\\
y=y\end{array}\right.\Longleftrightarrow \left\lbrace\begin{array}{l}x=\F{X}{y}\\y=y\end{array}\right.$, the jacobian of which is $\di\F{1}{y}$. We obtain
$$\int\int_{\O^*}(z'_x(x,y))^2\,dxdy=(\int_{\F{3}{4 }}^{1}y^3\,dy)(\int_0^{\F{1}{4}}\F{((\a'(X))^2}{(\a(X))^4}\,dX).$$ 
Since the function, $g$ defined by (\ref{2dg5}), is strictly increasing from $[1,\F{3}{2}]$ to 
$[0,\F{3}{4}]$, making a substitution, $\left\lbrace\begin{array}{l}X=g(u)\\
1\le u\le \a(\F{1}{4})\end{array}\right.\Longleftrightarrow \left\lbrace\begin{array}{l}u=\a(X)\\0\le X\le\F{1}{4}\end{array}\right.$ 
   yields
$$\int\int_{\O_*}(z'_x(x,y))^2\,dxdy=(\int_{\F{3}{4 }}^{1}y^3\,dy)( \int_1^{\a(\F{1}{4})}\F{(\a'(g(u)))^2g'(u)}{(\a(g(u)))^4}\,du).$$
Since $\a'(g(u))=\di\F{1}{g'(u)}$ and $\a(g(u))=u$, for $1<u<\a(\F{1}{4})$, we derive 
$$\int\int_{\O_*}(z'_x(x,y))^2\,dxdy=(\int_{\F{3}{4 }}^{1}y^3\,dy)(\int_1^{\a(\F{1}{4})}\F{1}{u^4g'(u)}\,du)\ge \F{175}{1024(\a(\F{1}{4}))^4}
\int_1^{\a(\F{1}{4})}\F{1}{g'(u)}\,du.$$
In view of $$\F{1}{g'(u)}=\F{\sqrt{(2u-1)(3-2u)}(\sqrt{(2u-1)(3-2u)}-8u^2+12 u-3)}{2(u-1)(-16u^3+32u^2-17u+3)}\ \mbox{$\sim\atop 1$}\ \ \F{1}{2(u-1)},$$
$$\mr{since}\ \int_1^{\a(\F{1}{4})}\F{1}{2(u-1)}\,du=+\iy,\ \mr{we\ obtain}\ \int\int_{\O^*}(z'_x(x,y))^2\,dxdy=+\iy.$$
Finally, the solution $z$ of the example 7 does not belong to $H^1(\O)$ and, therefore, the pro-\linebreak blem (\ref{pbh2}) is not well-posed. The reason why is probably that 
$(\u\,.\,\mb{n})(t)$, which vanishes in $t=\pi$, is not equivalent to $A(t-\pi)$, with $A\not=0$, in the neighborhood of $\pi$, that is to say it vanishes with an order greater than 1. On the contrary, the assumption $(\u\,.\,\tho_-)(t)$ negative is verified for $t=t_0$ and $t=\pi$.
\vspace{1cm}\\ 
 \textbf{References}\\[0.2cm]
\lk 1\rk\hs L. Ambrosio, Transport equations and Cauchy problems for
BV vector fields, Invent. Math., 158 (2004), pp. 227-260.\\[0.2cm]
\lk 2\rk\hs C. Bardos, Probl\`emes aux limites pour les \'equations
aux d\'eriv\'ees partielles du premier ordre \`a coefficients r\'eels;
Th\'eor\`emes d'approximations; Application \`a l'\'equations de
transport, Ann. Sci. \'Ecole Norm. Sup. (4), 3 (1970), pp. 185-223.\\[0.2cm]
 \lk 3\rk\hs J. M. Bernard, Steady transport equation in the case where the normal component of the velocity does not vanish on the boundary, SIAM J. Math. Anal., 
 Vol. 44, No. 2 (2012), pp. 993-1018.\\[0.2cm]
\lk 4\rk\hs F. Colombini and N. Lerner, Uniqueness of continuous
solutions for  vector fields, Duke Math. J., 111 (2002),
pp. 247-273.\\[0.2cm]
\lk 5\rk\hs R. J. DiPerna and P. L. Lions, Ordinary differential 
equations, transport theory and Sobolev spaces, Invent. Math., 98 
(1989), pp 511-547.\\[0.2cm]
\lk 6\rk\hs V. Girault and P.A. Raviart, \textit{Finite Element
    Approximation for Navier-Stokes Equations. Theory and Algorithms,}
  SMC 5, Springer-Verlag, Berlin, 1986.\\[0.2cm]
 \lk 7\rk\hs V. Girault and L.R. Scott, Analysis of two-dimensional
   grade-two fluid model with a tangential boundary condition,
   J. Math. Pures Appl., 78 (1999), pp. 981-1011.\\[0.2cm]
  \lk 8\rk\hs V. Girault and L.R. Scott, Finite-element discretizations of a two-dimensional grade-two fluid model, Mod\'el. Math. et Anal. Num\'er. 
  \tb{35} (2001), pp 1007-1053.\\[0.2cm]
  \lk 9\rk\hs V. Girault and L. Tartar, $L^p$ and $W^{1,p}$ regularity of the solution of a steady transport equation, C. R. Acad. Sci. Paris, Ser.1 348 (2010), 
  pp. 885-890.\\[0.2cm]
  \lk 10\rk\hs P. Grisvard, \textit{Elliptic Problems in Nonsmooth Domains},
  Pitman Monographs and Stu-\linebreak dies in Mathematics 24, Pitman, Boston, MA,
  1985.\\[0.2cm]
 \lk 11\rk\hs J.P. Puel and M.C. Roptin, Lemme de
  Friedrichs. Th\'eor\`eme de densit\'e r\'esultant du lemme de
  Friedrichs, Rapport de stage dirig\'e par C. Goulaouic, Dipl\^{o}me
  d'Etudes Approfondies, Universit\'e de Rennes, 1967.\\[0.2cm]
\lk 12\rk\hs N. J. Walkington, Convergence of the discontinuous
Galerkin method for discontinuous solutions, SIAM J. Numer. Anal., 42
(2005), pp. 1801-1817.
 \end{document}